\documentclass[journal]{IEEEtran}

\IEEEoverridecommandlockouts      

\usepackage{amsmath}
\allowdisplaybreaks[4]

\renewcommand{\baselinestretch}{1}

\usepackage{epsfig,tabularx,amsfonts,amssymb}
\usepackage{amsfonts}
\usepackage{epstopdf}
\usepackage{graphics} 
\usepackage{epic,eepic,epsfig}
\usepackage{mathrsfs}
\usepackage{subfigure}
\usepackage{calc,ifthen}
\usepackage{color}
\usepackage{bm}
\usepackage{float}
\usepackage{pstricks}
\usepackage[normalem]{ulem}
\usepackage{hyperref}
\usepackage[capitalize, nameinlink]{cleveref}
\newcommand{\stkout}[1]{\ifmmode\text{\sout{\ensuremath{#1}}}\else\sout{#1}\fi}

\usepackage{enumerate}

\usepackage{fancyhdr}

\usepackage{supertabular}

\newtheorem{definition}{Definition}
\newtheorem{problem}{\upshape OCP}
\newtheorem{lemma}{Lemma}

\newtheorem{remark}{Remark}
\newtheorem{assumption}{Assumption}

\newtheorem{proposition}{Proposition}

\newtheorem{corollary}{Corollary}

\def\lungo #1{\mathord{\buildrel{\lower3pt\hbox{$\scriptscriptstyle\frown$}}
\over #1 } }

\def\1D#1#2{{{\partial}\over{\partial #2}}#1}

\def\1d#1#2{{{d}\over{d #2}}#1}

\newcommand{{\R}}{{\mathbb{R}}}
\newcommand{{\C}}{{\mathbb{C}}}

\newcommand{\N}{{\mathbb{N}}}

\title{A Hierarchical Robust Control Strategy for Decentralized Signal-Free Intersection Management
}

\author{Xiao Pan, Boli Chen, Li Dai, Stelios Timotheou, and Simos A. Evangelou
\thanks{This work has been supported by the EPSRC Grant EP/N022262/1, and the National Natural Science Foundation of China under Grants 62173036 and 62122014, and partially funded by the European Union{'}s Horizon 2020 research and innovation programme under grant agreement No 739551 (KIOS CoE) and the Government of the Republic of Cyprus through the Directorate General for European Programmes, Coordination and Development.}
\thanks{X. Pan and S. A. Evangelou are with the Dept. of Electrical and Electronic Engineering, Imperial College London, London, UK {\tt\small (xiao.pan17@ic.ac.uk, s.evangelou@ic.ac.uk}).}
\thanks{B. Chen is with the Dept. of Electronic and Electrical Engineering, University College London, London, UK {\tt\small (boli.chen@ucl.ac.uk)}.}
\thanks{L. Dai is with the School of Automation, Beijing Institute
of Technology, Beijing, China {\tt\small (li.dai@bit.edu.cn)}.}
\thanks{S. Timotheou is with the Dept. of Electrical and Computer Engineering and the KIOS Research and Innovation Center of Excellence, University of Cyprus, Cyprus {\tt\small (timotheou.stelios@ucy.ac.cy)}.}}

\begin{document}

\thispagestyle{empty}
\setcounter{page}{0}
\begin{figure*}
\centering
\includegraphics[width=.9\textwidth]{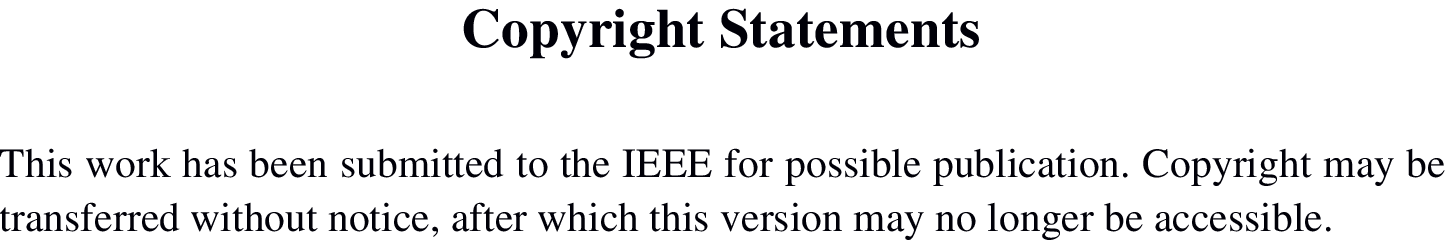}
\end{figure*}

\maketitle

\thispagestyle{fancy}
\chead{This work has been submitted to the IEEE for possible publication. Copyright may be transferred without notice, after which this version may no longer be accessible.}
\rhead{~\thepage~}
\renewcommand{\headrulewidth}{0pt}

\pagestyle{fancy}
\chead{This work has been submitted to the IEEE for possible publication. Copyright may be transferred without notice, after which this version may no longer be accessible.}
\rhead{~\thepage~}
\renewcommand{\headrulewidth}{0pt}

\begin{abstract}
The development of connected and automated vehicles is the key to improving urban mobility safety and efficiency. This paper focuses on cooperative vehicle management at a signal-free intersection with consideration of vehicle modeling uncertainties and sensor measurement disturbances. The problem is approached by a hierarchical robust control strategy in a decentralized traffic coordination framework where optimal control and tube-based robust model predictive control methods are designed to hierarchically solve the optimal crossing order and the velocity trajectories of a group of CAVs in terms of energy consumption and throughput. To capture the energy consumption of each vehicle, their powertrain system is modeled in line with an electric drive system. With a suitable relaxation and spatial modeling approach, the optimization problems in the proposed strategy can be formulated as convex second-order cone programs, which provide a unique and computationally efficient solution. A rigorous proof of the equivalence between the convexified and the original problems is also provided. Simulation results illustrate the effectiveness and robustness of the proposed strategy and reveal the impact of traffic density on the control solution. The study of the Pareto optimal solutions for the energy-time objective shows that a minor reduction in journey time can considerably reduce energy consumption, which emphasizes the necessity of optimizing their trade-off. Finally, the numerical comparisons carried out for different prediction horizons and sampling intervals provide insight into the control design.

\end{abstract}



{
\section*{Acronyms}
\noindent
\begin{tabular*}{0.49\textwidth}{@{}l@{\hspace{2mm}} @{\extracolsep{\fill}}l}
 CAV & Connected and Automated Vehicle\\
 CZ  & Control Zone\\
 FIFO & First-In-First-Out\\
 HRCS & Hierarchical Robust Control Strategy\\
 RMPC & Robust Model Predictive Control \\
 MZ & Merging Zone\\
 OCP & Optimal Control Problem\\
\end{tabular*}
}

\section{Introduction}
Rapid urbanization leads to increasing demand for mobility in transportation, and consequently leading to severe road congestion, higher energy consumption, and more traffic accidents. 
To address these problems, research studies have focused on cooperative vehicle management to improve traffic efficiency, reduce accidents and save energy. The development of connected and autonomous vehicle (CAV) and wireless communication technologies, such as vehicle-to-infrastructure and vehicle-to-vehicle communications, can improve traffic capacity under 
complex on-road traffic scenarios, such as signalized and unsignalized intersections~\cite{guanetti2018control,klein2016emergence,zhao2022enhanced}.
A comprehensive overview of signalized and signal-free intersection management strategies is presented in~\cite{zhong2020autonomous,namazi2019intelligent}, {in which the problem formulation technique can be mainly categorized into model-based~\cite{Chen2021TITS,hu2021TranC} and data-driven~\cite{Antonio2022TVT,Pozzi2020CDC,Bae2019ACC} approaches. The former relies on first-principle parametric models while the latter utilizes non-parametric learning approaches for control decision-making, but the majority of the data-driven methods focus on signalized intersections.}
This paper focuses on unsignalized intersection management {via a model-based technique} and deals with autonomous vehicle coordination through advanced vehicular communication systems (e.g., ehicle-to-infrastructure) to reduce traffic delays and energy usage.

Both centralized and decentralized control solutions are studied in the literature~\cite{vcakija2019autonomous,gholamhosseinian2022comprehensive}. Centralized approaches involve a central intersection controller that determines the optimal velocity profiles of all CAVs to safely cross the intersection aiming to satisfy objectives, such as safety and throughput maximization, and energy minimization~\cite{mihaly2020model,murgovski2015convex,fayazi2018mixed,castiglione2020cooperative,guney2020scheduling}. The trade-off between travel time and energy consumption minimization is investigated in~\cite{hadjigeorgiou2019optimizing}, where a hierarchical optimization approach is developed subject to First-In-First-Out (FIFO) constraints. 
\cite{riegger2016centralized} proposes a centralized model predictive control (MPC) method for a convexified modeling framework. 
The electric powertrain and various friction losses that were usually omitted are included in recent works, where \cite{hult2019optimal} involves turning maneuvers and \cite{chen2020optimal,liu2020high,pan2022convex} utilize convex modeling techniques to reduce the overall computational complexity. 
The utilization of such a realistic modeling framework with explicit powertrain and vehicle longitudinal dynamics has been shown to yield more accurate results (particularly with respect to energy consumption), as compared to other works that employ conventional lossless vehicle models. 
The work in \cite{pan2022convex} proposes a hierarchical optimization approach that reformulates the problem in a convex optimization framework
to jointly optimize both the travel time and powertrain energy consumption of battery electric vehicles. The upper-level of the hierarchical framework determines an optimized crossing order while the lower-level derives the optimal trajectories of all CAVs so that collision avoidance is guaranteed for the given order.  
In addition to the energy consumption and travel time that are commonly addressed, driving comfort is also considered in \cite{zhao2018multi}, where trade-offs between all these control objectives are investigated.

Although centralized schemes usually achieve optimal solutions, 
their implementation is challenging as 
they may suffer from high computational costs and lack resilience against the failure of the central controller. A more practical alternative is offered by decentralized control frameworks, where the velocity trajectory of individual CAVs is found locally by solving a local optimization problem associated with each vehicle~\cite{chalaki2022optimal,hult2018optimal,chalaki2021priority,hadjigeorgiou2022real}. An analytic optimization method is proposed in \cite{malikopoulos2018decentralized}, where the energy-optimal trajectory of each CAV is found individually by Pontryagin’s minimum principle. 
\cite{kumaravel2021optimal} presents a virtual platoon management rule, which forms the CAVs arriving at the intersection as a virtual platoon so as to enhance both the speed and throughput of CAVs.
In \cite{wu2019dcl}, the sequential movement of CAVs through the intersection is modeled using multi-agent Markov decision processes, while reinforcement learning is employed to find each velocity trajectory. Turning maneuvers are integrated into \cite{khoury2019practical,mirheli2019consensus}, indicating the capability of distributed coordinated frameworks in finding near-optimal solutions. 
More specifically, the work in \cite{khoury2019practical} solves the cooperative trajectory planning problem using vehicle-level mixed-integer non-linear programming, whereas the work in \cite{mirheli2019consensus} solves the problem by resorting to analytic optimization techniques.
More computationally efficient alternatives for decentralized approaches are heuristic control strategies \cite{zhang2015state,chouhan2018autonomous}, which however do not have optimality guarantees in most cases. 

Despite a rich literature on autonomous intersection management techniques, there are only a few works that deal with uncertainties entailed in the system. In \cite{chohan2019robust}, Kalman estimation techniques are utilized to cope with unreliable communication links (noise, packet drops, delays) in autonomous vehicle trajectory planning, and the coordination is accomplished by global optimization. A heuristic intersection management method has been proposed in~\cite{khayatian2018rim} from a robust control perspective to compensate for the effects of model mismatches and possible external disturbances. Recent work in \cite{vitale2022autonomous} considers vehicle position uncertainties resorting to the sensor and prediction errors in linear Gaussian motion models for a centralized coordination scheme.
To cope with the deviations in the trajectory following, a priority-aware resequencing mechanism is introduced in  \cite{chalaki2022priority} through indirect feedback in the system.
In \cite{chalaki2022robust}, a robust coordination scheme is proposed by incorporating date-driven and Gaussian process regression methods to learn the deviation from the nominal trajectories of CAVs.
However, the crossing orders in the works~\cite{chohan2019robust,khayatian2018rim} are predefined rather than optimized, and the powertrain dynamics and various friction losses are neglected in \cite{vitale2022autonomous,chalaki2022priority,chalaki2022robust}.
Based on some preliminary results presented in~\cite{pan2021decentralized}, this paper proposes a hierarchical robust control strategy (HRCS) for autonomous intersection crossing in a decentralized coordination framework, where the optimization problem is formulated in the {spatial domain} to avoid the free end-time problem and to utilize a convex optimization framework.
Moreover, both model and measurement uncertainties are considered in the state-space model and addressed by the developed hierarchical control framework.
This formulation enables the consideration of state-independent unmodeled longitudinal nonlinearities (e.g. gradient resistance caused by road slopes), and 
measurement noises of locations and velocities caused by sensors and environmental disturbances.
In particular, the hierarchical strategy involves a two-level optimization, in which the crossing order is optimized in an upper-level by solving an optimal control problem (OCP), which then guides a lower-level controller to derive optimal solutions in real-time.
To address the aforementioned uncertainties, a tube-based RMPC is designed for the lower-level coordination based on the robust invariant set centered along the nominal trajectory.
{In summary, this paper makes the following contributions:

\begin{itemize}
    \item This work introduces a novel robust
    and decentralized optimization-based autonomous intersection coordination framework, which takes into account vehicle modeling uncertainties and measurement noise. 
    \item 
    In addition to velocity trajectory optimization, the proposed scheme also finds the crossing order of the CAVs based on the local trajectory optimization and heuristic rules. This leads to significantly better solutions compared to the case when a predefined non-optimized (e.g., FIFO) crossing policy is enforced.
    \item 
    Inspired by the authors' prior work on centralized intersection coordination~\cite{pan2022convex}, computationally efficient solutions for the formulated decentralized coordinated scheme are derived by suitably relaxing the non-convex constraints and reformulating the associated optimization problems into convex second-order cone programs. A rigorous proof of the equivalence between the convexified and the original non-convex problem is also provided to complete the framework, which implies the robust satisfaction of the safety and operational constraints.
\end{itemize}
}

The rest of this paper is organized as follows. Section~\ref{sec:description} introduces a convex modeling framework of autonomous intersection crossing that includes the powertrain model of the CAV. The formulation of the decentralized HRCS is given in Section~\ref{sec:hierarchical}. Simulation results and discussion are presented in Section~\ref{sec:simulation}. Finally, concluding remarks are given in Section~\ref{sec:conclusions}.

\emph{Notation}: 
Let $\R$ and $\mathbb{N}$ denote the real and natural numbers, respectively. For $r_1\!\in\! \R$ and $n_1,n_2\!\in\!\mathbb{N}$, ${\R}_{\geq r_1}$, ${\R}_{> r_1}$, $\mathbb{N}_{>n_1}$, and $\mathbb{N}_{[n_1,n_2]}$ denote sets $\{r\!\in\!\R|r\!\geq\! r_1\}$, $\{r\!\in\!\R|r\!> \!r_1\}$, $\{n\in\mathbb{N}|n\!>\! n_1\}$, and $\{n\!\in\!\mathbb{N}|n_1\!\leq \!n\!\leq\!n_2\}$, respectively.
For $m\in \mathbb{N}_{>0} $, $X\in\R^{m}$ and $Y\in\N^{m}$ denote vectors with $m$ rows and all the elements being real and natural numbers, respectively.
For $m,n\in \mathbb{N}_{>0} $, $X\in\R^{m\times n}$ and $Y\in\N^{m\times n}$ denote matrices with $m$ rows and $n$ columns and all the elements being real and natural numbers, respectively.
Given sets $\mathbb{W},\,\mathbb{V}$,
Minkowski sum addition is noted as $\mathbb{W} \oplus \mathbb{V} = \{x + y|x \in \mathbb{W},y\in \mathbb{V}\}$, and Minkowski difference as $\mathbb{W} \ominus \mathbb{V} = \{x| x+V \subseteq \mathbb{W}\}$. 



\section{Problem Statement}
\label{sec:description}
\subsection{System Modeling and General Problem Setup}
\begin{figure}[t!]
\centering
\includegraphics[width=.9\columnwidth]{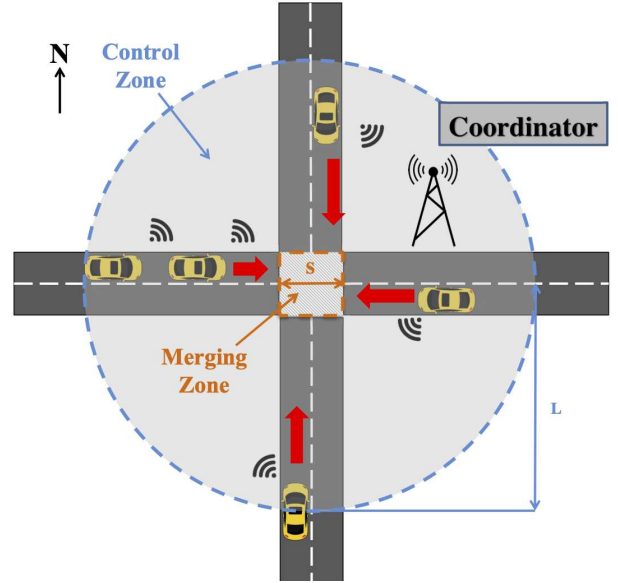}\\[-2ex]
\caption{The system architecture of autonomous intersection crossing problem.}
\label{fig:Intersection}
\end{figure}

The intersection coordination problem studied in the paper consists of a group of vehicles approaching a signal-free intersection, as shown in Fig.~\ref{fig:Intersection}, where lane changes and turning maneuvers are not allowed. All vehicles are assumed to be autonomous and connected, and there are no other non-autonomous road users (e.g. human-driven vehicles, cyclists, and pedestrians). As can be seen, the intersection is formed by two perpendicular roads with two lanes per road. Vehicles approaching the intersection will first enter a \textit{Control Zone} (CZ). The center of the intersection is called the \textit{Merging Zone} (MZ), where vehicles merge from different directions, and therefore lateral collision may occur. As it can be seen, the area of the MZ is considered as a square of side $S$ and the distance from the entry of the CZ to the entry of the MZ is $L$, with $L\!>\!S$. The intersection also has a coordinator that facilitates the exchange of information among the CAVs inside the CZ. Hence, in practice, $L$ is determined by the communication range capability of the CAVs and the coordinator. Note that the intersection coordinator is only used to streamline the communication process for the decentralized control scheme, and it is not involved in making control decisions. The decentralized control scheme will be specified in Section~\ref{sec:hierarchical}.
%
%

Let us denote by $N(t)\! \in\! \mathbb{N}_{>0}$ the total number of CAVs within the CZ at a given time $t\!\in\! \mathbb{R}_{>0}$ and the set $\mathcal{N}(t)\!\subseteq\!\mathbb{N}^{N(t)}$ to designate the crossing order in which the vehicles will enter the MZ. 
The determination of the crossing order will be elaborated in Section~\ref{section:upper_level}
The control target is to minimize the average energy consumption and travel time of $N(t)$ CAVs by finding in a decentralized manner the optimal sequence $\mathcal{N}(t)$ and the optimal speed trajectory for each vehicle from the entry of the CZ to the exit of the MZ. For notational convenience, we consider that $\mathcal{N}$ and $N$ are identical to $\mathcal{N}(t)$ and $N(t)$ in the rest of this paper, respectively.

\begin{definition}
Given an arbitrary CAV $h \in \mathcal{N}$, any CAV $i\in \mathcal{N},\,i\ne h$ can be categorized into one of the following subsets of $\mathcal{N}$ based on its physical location inside the CZ: 1) $\mathcal{C}_h$ collects vehicles traveling in the same direction as the $i$th vehicle; 2) $\mathcal{L}_h$ collects vehicles traveling in the perpendicular directions to the $i$th vehicle; 3) $\mathcal{O}_h$ collects vehicles traveling in the opposite direction to the $i$th vehicle.
\end{definition}

This research focuses on developing a coordination scheme in the {spatial domain}, which leads to the following two advantages: (a) the free end-time optimization problem in the time domain is avoided, and (b) the problem can be formulated using convex programs~\cite{pan2022convex}. Let $s$ denote the variable of traveled distance and $v_i(s)$ denote the velocity of the $i$th vehicle. The transformation from time to {spatial domain} is achieved by changing the independent variable $t$ to $s$ via $\frac{d}{ds}\!=\!\frac{1}{v}\frac{d}{dt}$. Thus, in the {spatial domain}, the distance traveled for each CAV is constant and equals to $L_i=(L+S+l_i),\,i\in \mathcal{N}$, and the travel time of each CAV can be readily obtained as a state variable, 
\begin{equation}\label{eq:dt}
  \frac{d}{ds}t_i =\frac{1}{v_i}  = \frac{1}{\sqrt{2E_i/m_i}}\,,
\end{equation}
where $l_i$ is the length of vehicle $i$, $E_i(s)\!=\!\frac{1}{2}m_i v_i^2(s)$ is the kinetic energy, and $m_i$ is the mass of the $i$th vehicle. As such, the required travel time of the $i$th CAV to cross the intersection is:
\begin{equation}\label{eq:J_ti}
J_{t,i}= t_i(L_i)-t_i(0)\,,
\end{equation}
where $t_i(0)$ is the arrival time of CAV $i$ at the CZ. Hereafter, $E_i(s)$ is introduced instead of $v_i$ for modeling the motion dynamics to cancel the nonlinearity due to the air drag. Considering $E_i(s)$ and vehicle longitudinal dynamics, the motion of vehicle $i$ can be described by,
\begin{equation}
\frac{d}{ds} E_i(s) = F_{w,i}(s) -F_{r,i}-\frac{2f_{d,i}}{m_i} E_{i}(s)\,, \label{eq:space_dv}   
\end{equation}
where $F_{w,i}(s)\!=\!F_{t,i}(s)\!+\!F_{b,i}(s)$ is the total force acting on the wheels with $F_{t,i}(s)$ and $F_{b,i}(s)$ being the powertrain driving force and the friction brake force acting on the wheels, respectively. $F_{r,i}\!=\!f_{r,i} {m_i} g$ is the rolling resistance force of CAV $i$ with coefficient $f_{r,i}$, and $f_{d,i}$ is the coefficient of air drag resistance of the $i$th CAV. It is further assumed that all vehicles are equipped with electric drives, and are subject to constraints:
\begin{align}
& \frac{g_{r,i}}{r_{w,i}}T_{\min,i} \leq F_{t,i}(s)\,, \leq \frac{g_{r,i}}{r_{w,i}}T_{\max,i} \label{eq:forwardforcebound}\\
& m_ia_{\min,i}-\frac{g_{r,i}}{r_{w,i}}T_{\min,i}\leq F_{b,i}(s) \leq 0\,, \label{eq:brakeforcebound}\\
& m_ia_{\min,i}\leq F_{w,i}(s) \leq \frac{g_{r,i}}{r_{w,i}}T_{\max,i}\,, \label{eq:wheelforcebound}
\end{align}
where $r_{w,i}$ is the wheel radius, $g_{r,i}$ is the fixed transmission gear ratio, $a_{\min,i}$ is the peak deceleration during emergency braking, and $T_{\min,i},\, T_{\max,i}$ are electric drive motor torque limits for generating and motoring, respectively. 

The energy cost of each CAV, ${J}_{b,i}$, is evaluated based on the tank-to-wheel energy path of the vehicle, which can be represented as a time integral of  ${P}_{b,i}(F_{t,i}, v_i)$, a function of vehicle force and speed (which correspond to motor torque and speed). As such, the energy usage in the {spatial domain} of each CAV is:
$
{J}_{b,i} = \int_{0}^{L_i} \frac{P_{b,i}(F_{t,i}, v_i)}{v_i(s)}ds.
$
By analogy to the common quadratic battery power model \cite{Han:aut19}, the following fitting model is chosen for $P_{b,i}$:
\begin{equation}\label{eq:Pbapprox}
\begin{aligned}
{P}_{b,i}(F_{t,i}, v_i)={b}_2F_{t,i}^2v_i+{b}_1F_{t,i}v_i+{b}_0v_i\,,
\end{aligned}
\end{equation}
where $b_0,\,b_1,\,b_2$ are fitting parameters obtained on the basis of a motor map. As such, the battery energy usage in the {spatial domain} of each CAV is:
\begin{equation}\label{eq:battery_usageapprox2}
J_{b,i} \!= \!  \int_{0}^{L_i}\!\!\!\!\!\left({b}_2F_{t,i}(s)^2\!+\!{b}_1F_{t,i}(s)\!+\!{b}_0\right)ds.
\end{equation}

For safety purposes, rear-end and lateral collision avoidance constraints and speed limits are required:
\begin{align}
& t_i(s)\!-\!t_h(s+l_h) \geq t_{\delta},\,\,\forall i\in\mathcal{C}_h,\label{eq:TTCapprox}\\
 & t_i(L)\geq t_h(L_h)\,,\,\,\, \forall i\in\mathcal{L}_h\,,\label{eq:lateral}\\
&\frac{1}{2}m_iv_{\min}^2 \leq \,E_i(s)\, \leq \frac{1}{2}m_iv_{\max}^2\,, \label{eq:Bound_E_i}
\end{align}
where $t_{\delta}$ is a minimum time gap enforced to prevent rear-end collision between vehicle $i$ and the vehicle immediately ahead, 
$v_{\min}$ is set to a sufficiently small positive constant to avoid singularity issues that would appear in \eqref{eq:dt} when $v_i{=}0$, and $v_{\max}$ is determined based on the infrastructure constraints and traffic regulations~\cite{hadjigeorgiou2022real}. 
Lateral collision constraint \eqref{eq:lateral} guarantees that the $i$th vehicle enters the MZ only after the $h$th vehicle leaves the MZ. 
For any CAV $i \in \mathcal{O}_h$, there is no interference between CAVs $h$ and $i$ inside the CZ. Hence, only the following constraint
\begin{equation}
t_i(L\!+\!S)\!>\!t_h(L\!+\!S),\, \,\,\forall i\in\mathcal{O}_h\,,
\label{eq:opposite}
\end{equation}
is required to fulfill the crossing order. 

The following assumptions are also needed to complete the modeling framework described above:
\begin{assumption}\label{ass:delay}
	All vehicle information (e.g., position, velocity, acceleration) can be measured through sensors, and the data can be transferred between each CAV and the coordinator without delays.
\end{assumption} 

\begin{assumption}\label{ass:initial}
For each CAV $i$, constraints \eqref{eq:TTCapprox}, \eqref{eq:Bound_E_i} and \eqref{eq:forwardforcebound} are inactive at $t_i(0)$.
\end{assumption}

Assumption~\ref{ass:delay} may not be valid for practical vehicular networks. On that occasion, it can be relaxed by using a worst-case analysis as long as the measurement and communication delays are bounded.
Assumption~\ref{ass:initial} is needed to ensure that all CAVs arriving at the CZ have feasible initial states and initial control inputs. 

Based on the above Assumptions~\ref{ass:delay}-\ref{ass:initial} and equations \eqref{eq:dt}-\eqref{eq:opposite}, the autonomous coordination problem can be formulated as an OCP. For an assumed crossing order $\mathcal{N}$ and 
the objective of minimizing the travel time \eqref{eq:J_ti} and energy consumption \eqref{eq:battery_usageapprox2}, the OCP is defined as
\begin{problem}\label{prob:ocp}
\begin{subequations}\label{eq:ocp}
\begin{align}
    & \mathop {\text{minimize}}\limits_{\mathbf{u}_{\text{nonlinear}}} \hspace{3mm}  W_1 \sum_{i=1}^{N}J_{t,i} +W_2 \sum_{i=1}^{N}J_{b,i} \label{eq:J_convex} \\
    &\textbf{s.t.}:  
    \eqref{eq:dt}, \eqref{eq:space_dv}, \eqref{eq:forwardforcebound}, \eqref{eq:brakeforcebound}, \eqref{eq:TTCapprox}, \eqref{eq:lateral}, \eqref{eq:Bound_E_i}, \eqref{eq:opposite}, 
   \vspace{-3mm}
\end{align}
\end{subequations}
\end{problem}
where 
\begin{multline*}
\mathbf{u}_{\text{nonlinear}}\!=\![F_{t,1}(s),F_{t,2}(s),\cdots,F_{t,N}(s), F_{b,1}(s),\\F_{b,2}(s),\cdots,F_{b,N}(s)]^{\top} \in \mathbb{R}^{2N},  
\end{multline*}
and $W_1, W_2 \in \R_{>0}$ are weighting factors. Note that the crossing order $\mathcal{N}$ is determined by another OCP that is introduced later in \cref{section:upper_level}.

\subsection{Convex Modeling Approach}
\cref{prob:ocp} is a non-convex optimization problem because of the non-affine equality dynamics \eqref{eq:dt}. Subsequently, we show that \cref{prob:ocp} can be convexified such that under certain assumptions the solution of the relaxed, convex problem is identical to the solution of the original non-convex \cref{prob:ocp}.

Towards this direction, the dynamics of $t_i$ are rewritten as 
\begin{align}
&  \frac{d}{d s}t_i(s) = \zeta_i(s)\,,\label{eq:syst4}\\
 & \zeta_i(s) \geq  \frac{1}{\sqrt{2E_i(s)/m_i}}, \label{eq:syst3}
\end{align}
which relax the original nonlinear differential equation into a linear differential equation and a convex constraint of the auxiliary control variable $\zeta_i(s)$.

The following Assumption~\ref{ass:fittingparameters}-\ref{ass:brakeforce} are needed to proceed with the proof of \cref{prop:tightness} defined later in this section.

\begin{assumption}\label{ass:fittingparameters}
The regression model \eqref{eq:Pbapprox} can find an accurate fitting of $P_{b,i}$ by $b_1$ and $b_2$ that comply with the condition
\begin{equation}\label{eq:fittingparameters}
b_1+2b_2F_{w,\min}>0,\,\forall i\in \mathcal{N},
\end{equation}
where 
$F_{w,\min}\!=\!\min\{m_i a_{\min,i}\}\!<\!0, \forall i\!\in\! \mathcal{N}$ representing the maximum available braking force of all CAVs.
\end{assumption}

\begin{figure}[t!]
\centering
\includegraphics[width=.465\columnwidth]{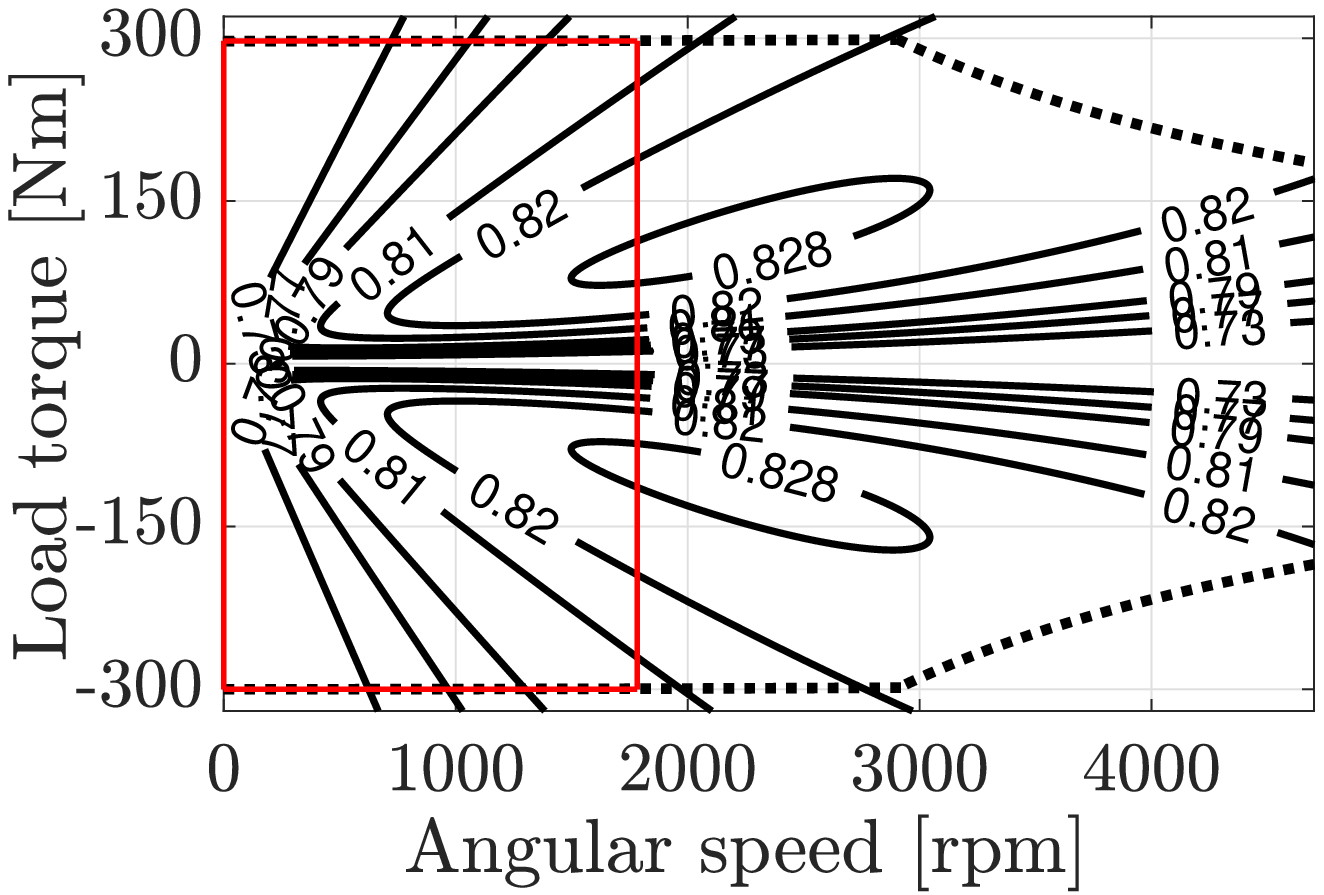}
\includegraphics[width=.515\columnwidth]{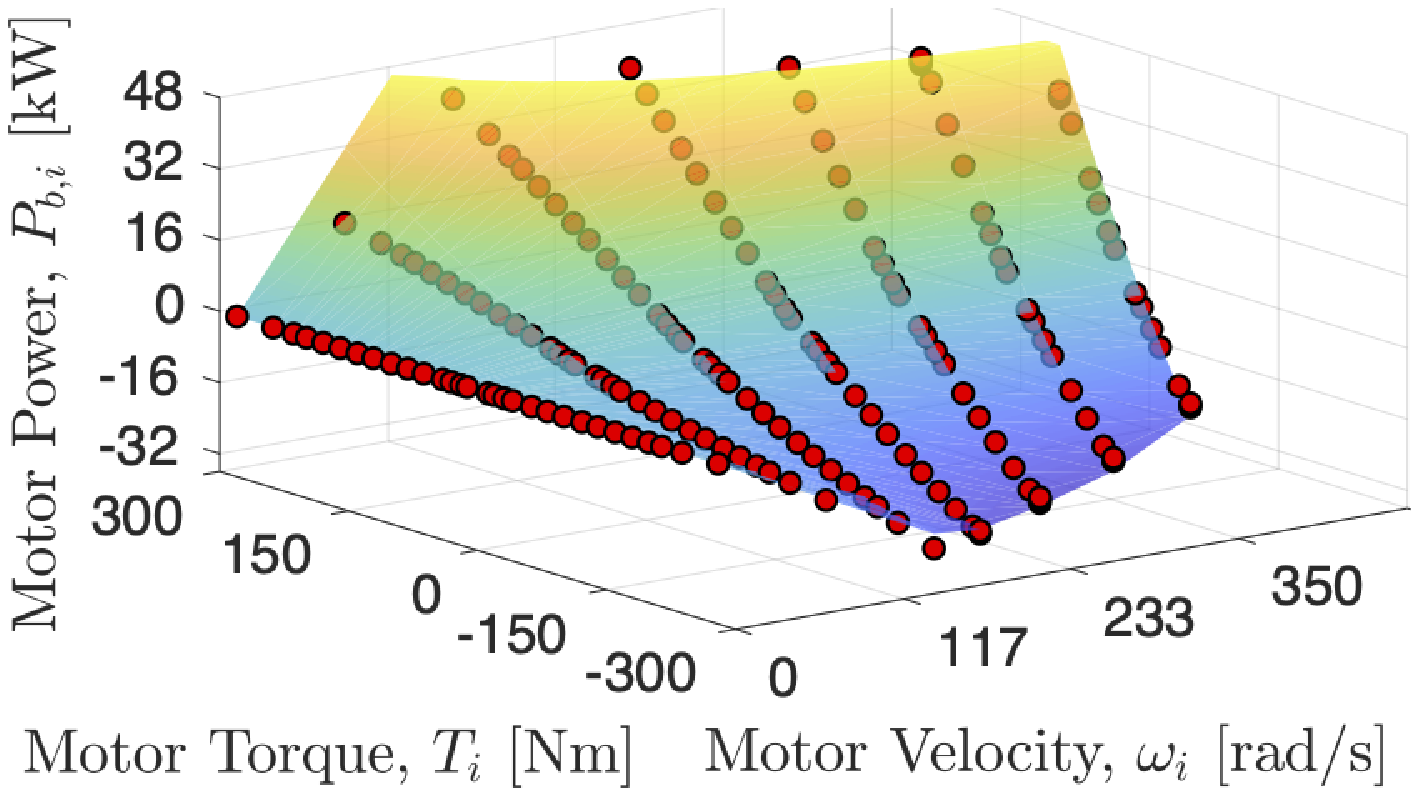}
\\[-2ex]
\caption{Left: efficiency map of the electric motor (positive torque indicates battery discharging and negative torque represents battery charging) with operational bounds (dotted lines). The area surrounded by red lines denotes the operational region for the feasible vehicle speed specified by \eqref{eq:Bound_E_i}. Right: nonlinear regression of the battery output power data (red dots, calculated based on the efficiency map shown in the left figure of Fig.~\ref{fig:motormap_Pbfitting}) by solving \eqref{eq:fit} with
an R-square fit of 99.25\% whereas the result of unconstrained fitting (without \eqref{eq:fit2}) is 99.53\%.
 }
\label{fig:motormap_Pbfitting}
\end{figure}

\cref{ass:fittingparameters} implies that $b_1$ has sufficiently large positive value compared to $b_2$. This
is consistent with the magnitude of the terms of the fitting model \eqref{eq:Pbapprox} based on their physical interpretation, where the fitting term involving $b_1$ is found to be the dominant one when the model \eqref{eq:Pbapprox} is fitted to data.
According to Assumption~\ref{ass:fittingparameters}, it is immediate to determine the fitting parameters $b_0,\,b_1,\,b_2$ by solving the following constrained optimization problem:
\begin{subequations}\label{eq:fit}
\begin{align}
    & \mathop {\min}\limits_{b_0,\,b_1,\,b_2 \in \mathbb{R}_{>0}} \hspace{3mm} ||P_{b,i}(F_{t,i}, v_i) - P^*_{b,i}||_2  \label{eq:fit1} \\
    &\textbf{s.t.}: \eqref{eq:fittingparameters}\,, \label{eq:fit2}
\end{align}
\end{subequations}
with $P^*_{b,i}$ the experimental battery output power data subject to certain motor torque and speed combinations. 
A representative example of the fitting \eqref{eq:Pbapprox} on the basis of the motor map given in \cite{chen:2019} is given Fig.~\ref{fig:motormap_Pbfitting}, which verifies the accuracy of the fitting model.  

\begin{assumption}\label{ass:brakeforce}
The friction brake $F_{b,i}$ is inactive, such that $F_{b,i}(s)=0,\,\forall i\in\mathcal{N},\,\forall s\in[0,L_i]$.
\end{assumption}
\cref{ass:brakeforce} holds in most cases by the fact that regenerative braking is naturally maximized in order to promote eco-driving, which targets minimum energy usage. The case when the assumption does not hold is analyzed in \cref{rem:activebrake} later in this section.

Under \cref{ass:brakeforce}, energy cost $J_{b,i}$ in \eqref{eq:battery_usageapprox2} is equivalent to
\begin{equation}\label{eq:battery_usageapprox3}
J_{w,i} = \int_{0}^{L_i}\left({b}_2F_{w,i}^2(s)\!+\!{b}_1F_{w,i}(s)\!+\!{b}_0\right)ds.
\end{equation}

Now, we have all ingredients to reformulate the non-convex problem \cref{prob:ocp} as a convex second-order cone program:
\begin{problem}\label{prob:ocp2}
\begin{subequations}\label{eq:ocp2}
\begin{align}
    & \mathop {\min}\limits_{\mathbf{u}} \hspace{3mm} W_1 \sum_{i=1}^{N} J_{t,i} +W_2 \sum_{i=1}^{N} J_{w,i} \label{eq:J_convex2} \\
    &\textbf{s.t.}: 
    \eqref{eq:space_dv}, \eqref{eq:wheelforcebound}, \eqref{eq:TTCapprox}, \eqref{eq:lateral}, \eqref{eq:Bound_E_i}, \eqref{eq:opposite}, 
    \eqref{eq:syst4},\eqref{eq:syst3}\,,
\end{align}
\end{subequations}
\end{problem}
where
$
\mathbf{u}\!=\![F_{w,1}(s),\cdots,F_{w,N}(s), \zeta_{1}(s),\cdots,\zeta_{N}(s)]^{\top} \!\in\! \mathbb{R}^{2N}  
$.

It is worth noting that the validity of the solution of \cref{prob:ocp2} relies on the tightness of \eqref{eq:syst3}, which is addressed by \cref{prop:tightness}.

\begin{proposition}\label{prop:tightness}
Under \cref{ass:fittingparameters} and \cref{ass:brakeforce}, the globally optimal solution of \cref{prob:ocp2} always finds the equality condition of \eqref{eq:syst3}, and therefore the solution of the relaxed, nonlinear convex problem \cref{prob:ocp2} is identical to the solution of the non-convex problem \cref{prob:ocp}.
\end{proposition}

\textit{Proof}. 
In view of \eqref{eq:J_convex2}, the multi-objective function of a single vehicle $i$ can be denoted by 
$
{J}_i=W_1J_{t,i}+W_2J_{w,i},
$
where 
$
J_{t,i}(\zeta_{i}(s)) = t_i(L_i)-t_i(0) = \int_0^{L_i}\zeta_i(s)$. Thus, $J_{t,i}$ only depends on $\zeta_{i}(s)$ whereas $J_{w,i}$ solely depends on $F_{w,i}(s)$, as it can be seen in \eqref{eq:battery_usageapprox3}.
Suppose that a feasible solution (denoted by the superscript $*$) of problem \cref{prob:ocp2}, with control signals $\zeta_i^*(s)$, $F_{w,i}^*(s)$, and the states $E_i^*(s)$, $t_i^*(s)$, is found for which it holds that $\zeta_i^*(s)\!>\!1/\sqrt{2(E_i^*(s))/m}$, and therefore we have
\begin{subequations}
\begin{align}
&\zeta_i^*(s)\!=\!\frac{1}{\sqrt{2(E_i^*(s))/m}}\!+\!\Delta \zeta_i(s),\, \text{with } \Delta \zeta_i(s)>0, \label{eq:zeta_star}\\ 
&t_i^*(s)\!=\!t_i(0)\!+\!\int_0^s \frac{1}{\sqrt{2(E_i^*(\xi))/m}}d\xi\!+\! \int_0^s \Delta \zeta_i(\xi) d\xi,
\end{align}
\end{subequations}
where the slack variable $\Delta \zeta_i$ inflates the travel time of the $i$th vehicle, and therefore it relaxes the rear-end and lateral collision avoidance constraints (e.g., if the velocity of CAV $i$ arriving at the CZ is much faster than that of CAV $i\!-\!1$ with $i,i\!-\!1\!\in\!\mathcal{N}$, it does not need to decelerate for CAV $i-1$, avoiding energy loss). 
It is also possible to construct an alternative feasible solution  $\left(\breve{\zeta}_i(s), \breve{E}(s), \breve{t}_i(s), \breve{F}_{w,i}(s)\right)$ with the same initial conditions ${E}_i^*(0)\!=\!\breve{E}_i(0)$ and ${t}_i^*(0)\!=\!\breve{t}_i(0)$, and
\begin{equation}\label{eq:breve_zeta}
    \breve{\zeta}_i(s)=\frac{1}{\sqrt{2(\breve{E}_i(s))/m_i}}=\zeta_i^*(s)\,,
\end{equation}
with $\breve{E}_i(s) \!<\! {E}_i^*(s)$, and therefore by integrating \eqref{eq:breve_zeta} $\breve{t}_{i}(s) \!=\!t_i^*(s),\forall s$. The alternative solution corresponds to the case when the $i$th vehicle slows down for the preceding vehicle without inflating the time variable (tightness of \eqref{eq:syst3} is guaranteed).
In this context, the kinetic energy difference between the two solutions can be obtained by 
\begin{equation}
\begin{aligned}
&\Delta E_i(s) =E_i^*(s)-\breve{E}_i(s)=  E_i^*(s)-\frac{1}{2}m_i\left(\frac{1}{\breve{\zeta}_i(s)}\right)^2 \\
&=E_i^*(s)-\frac{1}{2}m_i\frac{2E_i^*(s)}{(\sqrt{m_i}+\sqrt{2E_i^*(s)}\Delta \zeta_i(s))^2}> 0\,,
\end{aligned}
\end{equation}
where $\frac{1}{\breve{\zeta}_i(s)} \!=\! \frac{1}{{\zeta}_i^*(s)}$ owing to \eqref{eq:breve_zeta}, and therefore can be determined by \eqref{eq:zeta_star}.
Such a kinetic energy difference results by the deviation between $\breve{F}_{w,i}$ and ${F}^*_{w,i}$, $\Delta F_{w,i}(s)\! =\! {F}^*_{w,i}(s)\!- \!\breve{F}_{w,i}(s)$. The relationship between $\Delta E_{i}(s)$ and $\Delta F_{w,i}(s)$ can be found by \eqref{eq:space_dv}. 
For both solutions cases, by integrating both sides of the longitudinal dynamic equation, it holds that
\begin{subequations}\label{eq:xx12}
\begin{align}
&E^*_i(s) \!=\!E_i(0)\!+\!\!  \int_0^s\!\!\!\! F^*_{w,i}(\xi)d\xi \!- \!\!\!\int_0^s\!\! \!\!F_{r,i}d\xi \! +\! \varepsilon_i\! \!\int_0^s\! \!\!\!E_i^*(\xi)d\xi, \label{eq:xx1}\\
&\breve{E}_i(s) \!=\!E_i(0)\!+\! \! \int_0^s\!\!\!\! \breve{F}_{w,i}(\xi)d\xi \!-\!\! \int_0^s\!\!\!\! F_{r,i}d\xi \!+\! \varepsilon_i\!\!\int_0^s\!\! \!\!\breve{E}_i(\xi)d\xi, \label{eq:xx2}
\end{align}
\end{subequations}
where $\varepsilon_i \!=\!-\frac{2f_{d,i}}{m_i}$ is a constant.
By subtracting \eqref{eq:xx2} from \eqref{eq:xx1}, we obtain
\begin{equation}\label{eq:deltafsign}
\begin{aligned}
\Delta E_i(s)\!=\!\int_{0}^s \!\!\Delta F_{w,i}(\xi)d\xi +\varepsilon_i \!\!\int_{0}^s \!\!\Delta E_{i}(\xi)d\xi > 0 ,\,\forall s\,.
\end{aligned}
\end{equation}
Since $\varepsilon_i < 0$ and $\Delta E_i(s)>0$, then $\varepsilon_i \int_{0}^s \Delta E_{i}(\xi)d\xi <0$, which implies from \eqref{eq:deltafsign}
\begin{equation}\label{eq:deltaF}
\begin{aligned}
\int_{0}^s \Delta F_{w,i}(\xi)d\xi > 0\,.
\end{aligned}
\end{equation}

Let $J^*_i({F}_{w,i}^*,{\zeta}_i^*)$ and $\breve{J}_i(\breve{F}_{w,i},\breve{\zeta}_i)$ denote the cost for a single vehicle $i$ in both solution cases. Then, their difference can be calculated by 
\begin{equation}\label{eq:costdifference}
\begin{aligned}
J^*_i({F}_{w,i}^*,&{\zeta}_i^*) - \breve{J}_i(\breve{F}_{w,i},\breve{\zeta}_i)\\
&=  W_1J_{t,i}({\zeta}_i^*(s))+W_2J_{w,i}({F}_{w,i}^*(s)) \\ &-  W_1J_{t,i}(\breve{\zeta}_i(s))-W_2J_{w,i}(\breve{F}_{w,i}(s))\\
&=  W_2(J_{w,i}({F}_{w,i}^*(s))-J_{w,i}(\breve{F}_{w,i}(s)))\,.
\end{aligned}
\end{equation}
In virtue of the quadratic form of $J_{w,i}(\cdot)$, \eqref{eq:costdifference} can be rearranged as
\begin{equation}\label{eq:costdifference2}
\begin{aligned}
&J^*_i({F}_{w,i}^*,{\zeta}_i^*) - \breve{J}_i(\breve{F}_{w,i},\breve{\zeta}_i)\\
&=  W_2\int_0^{L_i}b_2{{F}_{w,i}^*}(s)^2+b_1{{F}_{w,i}^*}(s)+b_0\,ds\\
&-  W_2\int_0^{L_i}b_2{\breve{F}_{w,i}}(s)^2+b_1{\breve{F}_{w,i}}(s)+b_0\,ds \\
&=  W_2\int_0^{L_i}\left[b_2({F}_{w,i}^*(s)\!+\!\breve{F}_{w,i}(s))\!+\!b_1\right]\Delta F_{w,i}(s) ds\,.
\end{aligned}
\end{equation}
{ 
Let us define $g_i(s) \!=\! \frac{F^*_{w,i}(s)+\breve{F}_{w,i}(s)}{2}\in[F_{w,\min},F_{w,\max}]$ with $\,F_{w,\max}\!=\!\max\{\frac{g_{r,i}}{r_{w,i}}T_{\max,i}\}$, $\,F_{w,\min}\!=\!\min\{m_ia_{\min,i}\}$ $\forall i\in \mathcal{N}$ and $F_{w,\min}<0$ and $F_{w,\max}>0$, (27) can be rewritten as:
\begin{equation*}
\begin{aligned}
&J^*_i({F}_{w,i}^*,{\zeta}_i^*) - \breve{J}_i(\breve{F}_{w,i},\breve{\zeta}_i)\\
& = W_2b_1\int_0^{L_i}\Delta F_{w,i}(s) ds + W_2b_2\int_0^{L_i}2g_i(s)\Delta F_{w,i}(s) ds\\
& = W_2b_1\int_0^{L_i}\Delta F_{w,i}(s) ds + W_2b_2\int_0^{L_i}2g_i(s)\times \\
&\hspace{2cm}\left(\Delta F_{w,i}(s)-\Delta F_{w,\min}\right)+2g_i(s)\Delta F_{w,\min}ds
\end{aligned}
\end{equation*}
where $\Delta F_{w,\min}<0$ is the lower bound of $\Delta F_{w,i}(s)$.
Since the coefficient $b_2$ is positive, $g_i(s)$ is continuous, and the integrable function of $\left(\Delta F_{w,i}(s)-\Delta F_{w,\min}\right)\geq0,\forall s \in [0,L_i]$ does not change sign, it can be inferred from the Mean Value Theorem that there exists $\overline{F}_{w,i}\in[F_{w,\min},F_{w,\max}]$ such that
\begin{equation*}
\begin{aligned}
&J^*_i({F}_{w,i}^*,{\zeta}_i^*) - \breve{J}_i(\breve{F}_{w,i},\breve{\zeta}_i)\\
&=W_2(2b_2\overline{F}_{w,i}+b_1)\int_0^{L_i}\Delta F_{w,i}(s) ds-2W_2b_2\overline{F}_{w,i}\times\\
&\hspace{2mm} \Delta F_{w,\min}(L_i)+2W_2b_2\Delta F_{w,\min}\int_0^{L_i}g_i(s)ds
\end{aligned}
\end{equation*}
As $g_i(s) \in[F_{w,\min},F_{w,\max}]$, it follows:
\begin{equation}\label{eq:meanvalue}
\begin{aligned}
&J^*_i({F}_{w,i}^*,{\zeta}_i^*) - \breve{J}_i(\breve{F}_{w,i},\breve{\zeta}_i)\\
&\geq W_2(2b_2\overline{F}_{w,i}+b_1)\int_0^{L_i}\hspace{-1mm}\Delta F_{w,i}(s) ds-2W_2b_2F_{w,\min}\times\\
&\hspace{2mm} \Delta F_{w,\min}(L_i)+2W_2b_2\Delta F_{w,\min}\int_0^{L_i}F_{w,\min}ds\\
&=W_2(2b_2\overline{F}_{w,i}+b_1)\int_0^{L_i}\Delta F_{w,i}(s) ds
\end{aligned}
\end{equation}
}
From \eqref{eq:deltaF}, \cref{eq:meanvalue} implies that $J^*_i({F}_{w,i}^*,{\zeta}_i^*) > \breve{J}_i(\breve{F}_{w,i},\breve{\zeta}_i)$ if $2b_2\overline{F}_{w,i}+b_1>0$, which is guaranteed by \eqref{eq:fittingparameters} in \cref{ass:fittingparameters}: 
\begin{equation}
    \left(2b_2\overline{F}_{w,i}+b_1\right) > 2b_2F_{w,\min}+b_1 >0\,.
\end{equation}
Hence, given a solution set $\left(\zeta_i^*(s), F_{w,i}^*(s), E_i^*(s), t_i^*(s)\right)$ without holding the equality condition of \eqref{eq:syst3} (subject to a slack variable $\Delta \zeta_i(s)$), there always exists an alternative solution $\left(\breve{\zeta}_i(s), \breve{F}_{w,i}(s), \breve{E}_i(s), \breve{t}_i(s)\right)$ with guaranteed tightness that is more optimal in terms of the individual cost $J_i$, which in turn applies to all CAVs (i.e., $\forall i\!\in\! \mathcal{N}$). Thus, the proof ends. 

{ Note that Proposition 1 only valid if feasible solutions exist. To ensure the existence of the feasible solutions, let us introduce the following \cref{lemma:existance}.  }

{ 
\begin{lemma}
\label{lemma:existance}
There always exists a sufficiently small constant $\sigma\!>\!0$ and let $v_{\min}\!=\!\sigma$, then a feasible solution can be found. 
\end{lemma}

\textit{Proof.} Consider a candidate solution $\left(\zeta_i^*(s),F_{w,i}^*(s),E_i^*(s),t_i^*(s)\right)$ that violates the lateral collision avoidance constraint \eqref{eq:lateral}, then we have,
\begin{equation*}
\begin{aligned}
    t_i^*(L) &\!<\! t_h(L_h)\\
  \Rightarrow  t_i(0)\!+\!\int_0^{L}\zeta_i^*(s)ds&\!<\!t_h(L_h)\\
  \Rightarrow  \int_0^{L}\frac{1}{\sqrt{2E^*_i(s)/m_i}}ds &\!<\!t_h(L_h)\!-\!t_i(0)
\end{aligned}
\end{equation*}
where $t_h(L_h)$ and $t_i(0)$ are preknown constants to vehicle $i$ in the decentralized framework.
Since ${E}_{i}$ is bounded by \eqref{eq:Bound_E_i} and we can conclude from \cref{lemma:existance} that
the vehicle $i$ can maintain the speed at $\sigma$, such that
as $\sigma\!\rightarrow\!0$ we have
${1}/{\sqrt{2E^*_i(s)/m_i}}\!=\!{1}/{\sigma}\rightarrow \infty$. As such, the lateral collision avoidance constraint can be satisfied as follows:
\begin{equation*}
\begin{aligned}
\int_0^{L}\!\!\!\frac{1}{\sqrt{2E^*_i(s)/m_i}}ds \!=\!\int_0^{L}\frac{1}{\sigma}ds &\geq t_h(L_h)\!-\!t_i(0).
\end{aligned}
\end{equation*}

Similarly, the satisfaction of the rear-end collision avoidance constraint \eqref{eq:TTCapprox} and the constraint to fulfill the crossing order \eqref{eq:opposite} can be justified. 

From Lemma~\ref{lemma:existance}, by solving the convex optimization problem with
continually reduced lower bound of the speed, there always exists a $\sigma$ and a corresponding feasible solution with ${E}^*_i(s)\!=\!\frac{1}{2}m\sigma^2$ that allows vehicle $i$ to operate at an almost static condition such that collision avoidance constraints is fulfilled, and therefore the equality condition \eqref{eq:syst3} holds.

}

\begin{remark} [Active friction brakes]
Solving \cref{prob:ocp2} may yield a solution trajectory where $F_{w,i}(s) \!<\! \frac{g_{r,i}}{r_{w,i}}T_{\min,i}(s)$ for some $s$ (regenerated braking power is saturated). In this circumstance, friction brake is invoked to meet the total force demand at the wheels, $F_{w,i}(s) \!=\! \frac{g_{r,i}}{r_{w,i}}T_{\min,i}(s) \!+\!F_{b,i}(s)$. As such, the equivalence between \cref{prob:ocp} and \cref{prob:ocp2} (see \cref{prop:tightness}) is no longer guaranteed due to the discrepancy between the energy costs in both OCPs (regeneration of friction brakes is assumed in \cref{prob:ocp2}), and therefore, the optimality of \cref{prob:ocp2} may be compromised in such a case. Nevertheless, the equality of \eqref{eq:syst3} holds invariably (as inferred from the above proof), which, in turn, ensures the feasibility of the convex optimization solution.
\label{rem:activebrake}
\end{remark}
%

\section{Hierarchical Robust Control Strategy}
\label{sec:hierarchical}
\begin{figure}[t!]
\centering
\includegraphics[width=1\columnwidth]{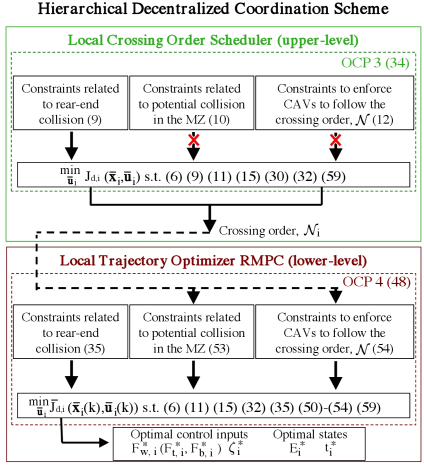}\\[-2ex]
\caption{{ The schematic of the decentralized HRCS. The objective functions in \cref{prob:ocp3} and \cref{prob:ocp4} take $J_{d,i}$ \eqref{eq:J_decentralized} and $\bar{J}_{d,i}$ \eqref{eq:decentralized_objective}, respectively.}}
\label{fig:Hierarchical_Scheme}
\end{figure}

This section introduces the novel decentralized control strategy (i.e., HRCS), where the convex problem formulated in~\cref{prob:ocp2} plays a key role in its formulation. As shown in Fig.~\ref{fig:Hierarchical_Scheme}, the scheme is composed of an upper-level crossing order scheduler and a lower-level trajectory optimizer RMPC, deployed in a hierarchical and decentralized manner. 
{ Every time a new vehicle $i$ arrives at the entry point of the CZ, its local controller determines the optimal passing order $\mathcal{N}_i$ by solving an OCP~\ref{prob:ocp3} as will be described in~\cref{section:upper_level}. The tube-based RMPC embedded in each CAV finds its individual optimal trajectory to cross the signal-free intersection.
%
Note that the intersection coordinator works as an information relay without making any control decisions. Specifically, the coordinator will collect the new crossing order $\mathcal{N}_i$ and transmit it all vehicles inside the CZ such that their lower-level controller will update the collision avoidance constraints in the lower-level optimization problem based on an updated crossing order.}

Prior to the introduction of the individual algorithms in the two layers, let us first introduce some preliminaries that are used across the two layers. 
Considering a sampling distance interval $\Delta s\!\in\! \mathbb{R}_{>0}$, and without loss of generality it is assumed that $L_i\!=\!\alpha\Delta s$, $L\!=\!\alpha_1\Delta s,\,S\!=\!\alpha_2\Delta s,\,l_i\!=\!\alpha_3\Delta s,\,\alpha,\alpha_1,\alpha_2,\alpha_3\!\in\!\mathbb{N}_{>0}, \alpha\!=\!\alpha_1\!+\!\alpha_2+\!\alpha_3,\alpha_1\!>\!\alpha_2\!>\!\alpha_3$. 

The objective function in a discretized form for vehicle $i$ is formulated based on \eqref{eq:J_convex2} as
\begin{equation}\label{eq:J_decentralized}
{J}_{d,i}= \sum_{k=0}^{\alpha-1}[W_1(b_2F_{w,i}^2(k)+b_1F_{w,i}(k)+b_0)
+W_2\zeta_{i}(k)]\Delta s,
\end{equation}

The discretized system equations \eqref{eq:space_dv} and \eqref{eq:syst4}, and the related constraints can be rewritten as:
\begin{subequations}\label{eq:discrete_dynamics}
\begin{align}
&x_{i}(k+1)= A_ix_{i}(k)\!+\!B_iu_{i}(k)+\omega_{i}(k)\,, \label{eq:discrete_dynamics_state}\\
&y_{i}(k)=C_ix_{i}(k)\!+\!\nu_{i}(k)\,, \label{eq:discrete_dynamics_output}\\
& x_{i} = \begin{bmatrix}
E_{i}\\
t_{i}
\end{bmatrix}\!,\quad 
u_{i} = \begin{bmatrix}
F_{w,{i}}-F_{r,i}\\
\zeta_i
\end{bmatrix}\!, \quad
A_i = \begin{bmatrix}
\mathrm{e}^{\frac{-2f_{d,i}\Delta s}{m_i}} & 0\\
0&1
\end{bmatrix}\!, \nonumber \\
& B_i\!=\! \begin{bmatrix}
-\frac{m_i(\mathrm{e}^{-\frac{2f_{d,i}\Delta s}{m_i}}-1)}{2f_{d,i}} & 0\\
0 & \Delta s
\end{bmatrix}\!, \quad
C_i \!=\! \begin{bmatrix}
1 & 0\\
0 & 1
\end{bmatrix}\!,
\nonumber
\end{align}
\end{subequations}
where $k\!=\!0,1,\ldots,\alpha\!-\!1$, $\omega_i(k)\!=\![\omega_{E,i}(k),\,0]^{\top}$ with a closed convex set bounding on process noise $\mathbb{W}_{i}\!=\!\{\omega_{E,i}(k)\!\in\!\mathbb{R}\!:\!||\omega_{E,i}(k)||\!\leq\!\overline{\omega}_{E,i}\}$ and $\overline{\omega}_{E,i}\!\in\!\mathbb{R}_{>0}$, $\nu_i(k)\!=\![\nu_{E,i}(k),\,\nu_{t,i}(k)]^{\top}$ with a closed convex set bounding on measurement noise $\mathbb{V}_{i}\!=\!\{\nu_i(k)\!\in\!\mathbb{R}^2\!:\!||\nu_{i}(k)||\!\leq\!\overline{\nu}_i\}$ and $\overline{\nu}_i\!\in\!\mathbb{R}_{>0}$, and $y_{i}$ is the measured signal of the $i$th vehicle's status, obtained by onboard sensors. The consideration of $\omega_{i}$ and $\nu_{i}$ involves, respectively, the state-independent unmodeled vehicle longitudinal dynamics (e.g. gradient resistance caused by road slopes), and the measurement noises from the sensors and environment disturbances. Moreover, the nominal system of the actual dynamics system \eqref{eq:discrete_dynamics_state} embedded in the designed RMPC controller can be written as:
\begin{equation}\label{eq:discrete_nominal_dynamics}
\bar{x}_{i}(k+1)= A_i\bar{x}_{i}(k)\!+\!B_i\bar{u}_{i}(k)\,,
\end{equation}
where the disturbance term $\omega_i$ 
is neglected, $\bar{x}_i(k)$ is the nominal state, and $\bar{u}_i(k)$ is the nominal control input.

In the decentralized framework, the coordinator assigns a unique identity to each CAV when the vehicle enters the CZ of the intersection based on the vehicle's arrival time and its entering direction. The identity that the coordinator assigns to each vehicle that arrives at the CZ is defined as $(i,d_i,\mathcal{I}_{i}(k))$, where $i\!\in\!\mathcal{N}_i$ is the optimal crossing order of vehicle $i$ to enter the MZ as will be described in OCP~\ref{prob:ocp3}, $d_i$ is an index denoting the traveling direction (north, south, east and west), and
$\mathcal{I}_{i}(k)$ is the information set generated by the coordinator at step $k$ by collecting the past optimal nominal state sequence of vehicle $i$:
\begin{equation}
\begin{aligned}
\mathcal{I}_{i}(k) &= [\mathbf{\bar{x}}_{i}(0),\mathbf{\bar{x}}_{i}(1),\cdots,\mathbf{\bar{x}}_{i}(k)] \in \mathbb{R}^{(N_p+1)\times 2k},\,
\end{aligned} 
\end{equation}
where $\mathbf{\bar{x}}_{i}(k)\!=\![\mathbf{\bar{E}}_{i}(k),\mathbf{\bar{t}}_{i}(k)]\!\in\!\mathbb{R}^{(N_p+1)\times 2}$ is the optimal nominal state sequence computed by the local controller on vehicle $i$ at step $k$ within a predefined receding horizon $N_p$, and 
\begin{align*}
\mathbf{\bar{E}}_{i}(k)&\!=\!\left[\bar{E}_{i}(k|k),\bar{E}_{i}(k+1|k),\ldots,\bar{E}_{i}(k+N_p|k)\right]^{\top}\!, \\
\mathbf{\bar{t}}_{i}(k)&\!=\!\left[\bar{t}_{i}(k|k),\bar{t}_{i}(k+1|k),\ldots,\bar{t}_{i}(k+N_p|k)\right]^{\top}\!.
\end{align*}

\subsection{Crossing Order Scheduler}
\label{section:upper_level}
Considering the single-lane road depicted in Fig.~\ref{fig:Intersection}, CAVs approaching the intersection from the same direction will enter and leave the MZ in the same order they arrive at the CZ. Therefore, the major challenge stems from the collision avoidance and prioritizing constraints for CAVs merging from different directions, i.e., \eqref{eq:lateral} and \eqref{eq:opposite}, which depend on the crossing order. To determine an optimal crossing order without invoking exhaustive search, we define a virtual coordination problem in the discrete-time form:

\begin{problem}\label{prob:ocp3}
\begin{subequations}\label{eq:upper_OCP}
\begin{align}
& \hspace{2mm}\mathop {\min}\limits_{\mathbf{\bar{u}_i}} \hspace{3mm} {J}_{d,i}(\bar{\mathbf{x}}_{i},\bar{\mathbf{u}}_{i}) \label{eq:upper_OCP_J} \\
    \textbf{s.t.}:\,\, &\bar{x}_i(0)={x}_{i}(0),\\
    &\bar{x}_{i}(k+1)= A_i\bar{x}_{i}(k)\!+\!B_i\bar{u}_{i}(k), \label{eq:nominal_dynamics} \\
    & \eqref{eq:space_dv}, \eqref{eq:wheelforcebound}, \eqref{eq:TTCapprox}, \eqref{eq:Bound_E_i}, 
    \eqref{eq:syst4},\eqref{eq:syst3}\,,\,k=0,1,\ldots,\alpha-1\,, \nonumber
\end{align}
\end{subequations}
\end{problem}
where $\mathbf{\bar{x}}_i\!=\![\bar{x}_i(0),\bar{x}_i(1),\ldots,\bar{x}_i(\alpha)]^{\top}\!\in\!\mathbb{R}^{(\alpha+1)\times 2}$, $\mathbf{\bar{u}}_i\!=\![\bar{u}_i(0),\bar{u}_i(1),\ldots,\bar{u}_i(\alpha-1)]^{\top}\!\in\!\mathbb{R}^{\alpha\times 2}$, \eqref{eq:nominal_dynamics} collects the nominal dynamics \eqref{eq:discrete_nominal_dynamics}, and the initial value of the nominal state $\bar{x}_i(0)$ takes the value of
the actual initial condition $x_i(0)$ since the disturbances are ignored.

{
Each vehicle needs to find the ideal (non-conservative) MZ entry and exit times, $t_{i,in}$ and $t_{i,out}$ by solving OCP~\ref{prob:ocp3} at the entry point of the CZ. Without considering the lateral collision avoidance constraints. The main idea is that for those vehicles that have potential lateral collisions, their crossing orders are based on their entry times at the MZ, while if there is no interference between the CAVs, their orders are based on their exit times at the MZ. As such, the implementation steps are performed as follows.

\begin{enumerate}[Step 1:]
\item Given $t_{i,in}$ and $t_{i,out}$ for vehicle $i \!\in\! \mathcal{A}_i$ with $\mathcal{A}_i=\{1,2,3\cdots,i\}$ the arrival order at the CZ, and a previously designated crossing order $\mathcal{N}_{i-1} \in\mathbb{N}^{(i-1)}$ {(i.e., $\mathcal{N}_{i-1}\!=\!\mathcal{N}(t),\forall t\!\in\!\!\![t_{i-1}(0),t_{i}(0))$)}, which is a permutation of $\mathcal{A}_{i-1}$. For the sake of further discussion, let us consider ${N}_{i-1}^{k}$ the $k$th element in $\mathcal{N}_{i-1}$.
\item 
The $i$th arriving vehicle receives the entry position and intention of the last vehicle, ${N}_{i-1}^{i-1}$ in $\mathcal{N}_{i-1}$ from the intersection coordinator, and then determines if there is lateral collision potential or not. 

\item (i) If there exists ``lateral collision potential", the new crossing order $\mathcal{N}_i$ is determined at vehicle $i$ by evaluating the ideal entry time of the $i$th arriving vehicle and the vehicle ${N}_{i-1}^{i-1}\in\mathcal{N}_{i-1}$ at the MZ:
\begin{equation*}
\mathcal{N}_i\!=\!
\left\{\!\!\!\!
\begin{array}{lll} 
\{{N}_{i-1}^{1},\ldots,{N}_{i-1}^{i-2},{N}_{i-1}^{i-1},i\}\!,\,t_{i,in}\!\geq\! t_{{N}_{i-1}^{i-1},in} ,\\
\{{N}_{i-1}^{1},\ldots,{N}_{i-1}^{i-2},i,{N}_{i-1}^{i-1}\}\!,\,t_{i,in}\!<\! t_{{N}_{i-1}^{i-1},in} ,
\end{array}\right.
\nonumber \\    
\end{equation*}
(ii) If there exists ``no lateral collision potential", the new crossing order $\mathcal{N}_i$ is determined at vehicle $i$ by evaluating the exit time of the $i$th arriving vehicle and the vehicle ${N}_{i-1}^{i-1}\in\mathcal{N}_{i-1}$ at the MZ:
\hspace{-5mm}
\begin{equation*}
\begin{aligned}
\mathcal{N}_i\!=\!  \left\{ \!\!\!\! \begin{array}{lll} 
\{{N}_{i-1}^{1},\ldots,{N}_{i-1}^{i-2},{N}_{i-1}^{i-1},i\}\!,\,t_{i,out}\!\geq \!t_{{N}_{i-1}^{i-1},out} ,\\
\{{N}_{i-1}^{1},\ldots,{N}_{i-1}^{i-2},i,{N}_{i-1}^{i-1}\}\!,\,t_{i,out}\!<\! t_{{N}_{i-1}^{i-1},out} ,
\end{array}\right.
\nonumber \\  
\end{aligned}
\end{equation*}
\item The coordinator receives $\mathcal{N}_i$ from vehicle $i$ and if $\mathcal{N}_i \ne \{\mathcal{N}_{i-1},i\}$, $\mathcal{N}_i$ is transmitted to vehicle ${N}_{i-1}^{i-1}\in\mathcal{N}_{i-1}$ and $(i+1)$th arriving vehicle. Otherwise, $\mathcal{N}_i$ is transmitted to $(i+1)$th arriving vehicle only.
\end{enumerate}

Given the crossing order, $\mathcal{N}_i$, the lateral collision avoidance constraints in the lower-level optimization are set up correspondingly to ensure safety.

}

\subsection{Trajectory Optimizer}
\label{section:lower_level}
Before the description of the lower-level controller, some preliminaries that are utilized to construct the collision constraints during the receding horizon windows are introduced first.
Based on the identity $(i,d_i,\mathcal{I}_{i}(k))$, the coordinator exchanges the information set $\mathcal{I}_{i}(k)$ with the associated vehicles (see below) to enable collision avoidance constraints to be established for each local controller according to \eqref{eq:TTCapprox} \eqref{eq:lateral} and \eqref{eq:opposite}. Depending on the crossing order of the CAV $i$, the information it requires from the coordinator to enable local control  is defined as follows:
\begin{enumerate}
    \item no information is required, if $i=1$,
    \item $\mathcal{I}_{i-1}(k_{i-1})$, if ${i}\in \mathcal{C}_{i-1}$,
    \item $\mathcal{I}_{i-1}(k_{i-1})$ and $\mathcal{I}_{h}(k_{h})$ with $h<i$, where $h$ stands for the CAV immediately ahead of CAV $i$, if $i \notin \mathcal{C}_{i-1}$ and $i \in \mathcal{C}_{h}$,
\end{enumerate}
where $k_{i-1}$ and $k_h$ are the corresponding distance step of CAV $i-1$ and $h$, respectively, when the $i$th CAV is at step $k$, that are ${t}_{h}(k_h \Delta s)\!=\!{t}_{i-1}(k_{i-1} \Delta s)\!=\!{t}_{i}(k \Delta s)$. If $i\!\in\! \mathcal{C}_{h}$, the rear-end collision constraint \eqref{eq:TTCapprox} is reformatted for the RMPC-based decentralized HRCS as follows:
\begin{equation}
{t}_{i}(k+j+1|k)-\bar{t}_{h}((k-\alpha_3)+j+1|k-\alpha_3) 
>t_{\delta} ,
\label{eq:decentralized_TTCapprox}    
\end{equation}
where $j=0,1,\cdots,N_p-1$, and $\bar{t}_{h}((k-\alpha_3)+j+1|k-\alpha_3)$ is the historical information stored in $\mathcal{I}_{h}(k_h)$. In addition to the rear-end collision avoidance constraint, the constraints concerning the lateral collision avoidance and prioritization for CAV $i\in\{\mathcal{L}_h,\mathcal{O}_h\}$, i.e., \eqref{eq:lateral} and \eqref{eq:opposite}, will be reformatted and described latter in \eqref{eq:decentralized_objective}-\eqref{eq:opposite_decentralize_format} in this section.

After obtaining the optimal crossing order $\mathcal{N}_i$ from OCP~\ref{prob:ocp3}, the trajectory of CAV $i$ can be optimized by solving an RMPC problem that replaces the objective function \eqref{eq:upper_OCP_J} 
with the later defined objective function \eqref{eq:decentralized_objective}-\eqref{eq:si_opposite} to satisfy conditions \eqref{eq:lateral} and \eqref{eq:opposite}, and to cope with the additive disturbance $\omega_i$ and $\nu_i$ in the actual system \eqref{eq:discrete_dynamics}. This work adopts an RMPC with output feedback and robust invariant tubes based on the nominal system \eqref{eq:discrete_nominal_dynamics}. An observer is firstly designed to measure the state of each vehicle. Then, we design robust invariant tubes based on the nominal system \eqref{eq:discrete_nominal_dynamics} to bound the nominal state and input. As such, by applying the optimal nominal input solved in each step $k$, the actual vehicle dynamics constraints $\{x_i\!\in\!\mathbb{X}_{i},u_i\!\in\!\mathbb{U}_{i},(x_i,u_i)\!\in\!\mathbb{X}_{i}\!\times\! \mathbb{U}_{i}\}$ can always be satisfied even with the effects of the admissible disturbance sequences $w_{i}$ and $\nu_{i}$. The set $\mathbb{X}_{i}$ collects constraints \eqref{eq:Bound_E_i} \eqref{eq:decentralized_TTCapprox} {as well as lateral collision avoidance and prioritization constraints \eqref{eq:lateral_decentralize_format}-\eqref{eq:opposite_decentralize_format} as described latter in this section}, and the set $\mathbb{U}_{i}$ collects \eqref{eq:wheelforcebound}, and the set $\mathbb{X}_{i}\times \mathbb{U}_{i}$ collects mixed state and input constraint \eqref{eq:syst3}.

The states of the system \eqref{eq:discrete_dynamics_state} are estimated by a Luenberger observer with dynamics described by:
\begin{equation}\label{eq:observer_dynamics}
\begin{aligned}
\hat{x}_i(k+1) &= A_i\hat{x}_i(k)+B_iu_i(k)+L_{c,i}(y_i(k)-\hat{y}_i(k))\,,\\
\hat{y}_i(k)&=C_i\hat{x}_i(k)\,,
\end{aligned}
\end{equation}
where $\hat{x}_{i}$ and $\hat{y}_{i}$ are the estimations of state and output, respectively, and $L_{c,i}$ is the observer gain such that the eigenvalues of $A_{L,i}(=\! A_i\!-\!L_{c,i}C_i)$ satisfy the condition $\lambda(A_{L_i})\!<\!1$.
Considering the estimation error
\begin{equation}\label{eq:estimation_error}
\tilde{x}_i(k) =x_i(k)-\hat{x}_i(k)\,,
\end{equation}
it can be readily shown that its dynamics are governed by 
\begin{equation}\label{eq:state_estimation_error}
\tilde{x}_i(k+1) = A_{L,i}\tilde{x}_i(k)+\tilde{\delta}_i(k)\,,
\end{equation}
where $\tilde{\delta}_i(k)\!=\!\omega_i(k)\!-\!L_{c,i}\nu_i(k)$. From the bounds of the disturbances, it follows that 
\begin{equation}
\tilde{\delta}_i(k)\in \tilde{\Delta}_{i}=\mathbb{W}_{i}\ominus L_{c,i}\mathbb{V}_{i}\,.
\end{equation}

Since $A_{L,i}$ is stable, a robust invariant tube $\tilde{\mathbb{S}}_{i}$ can be found such that if $\tilde{x}_i(0)\in \tilde{\mathbb{S}}_{i}$, $\tilde{x}_i(k)$ remain in the tube for all $\tilde{\delta}_i(k)$.
In view of \eqref{eq:state_estimation_error}, it holds that $A_{L,i}\tilde{\mathbb{S}}_{i}\oplus \tilde{\Delta}_{i} \subseteq\tilde{\mathbb{S}}_{i}$ and for the observer gain $L_{c,i}$, the robust invariant tube $\tilde{\mathbb{S}}_{i}$ that includes the effects of $\tilde{\delta}_i(k)\in \tilde{\Delta}_{i}$ can be computed as:
\[
\tilde{\mathbb{S}}_{i}=(1-\mu_{L,i})^{-1}\bigoplus_{j=0}^{n_{L,i}-1}A_{L,i}^{j}\tilde{\Delta}_{i}\,,
\]
where $n_{L,i}\!\in\!\mathbb{N}_{>0}$ is a finite integer and $\mu_{L,i}\!\in\![0,1)$ such that $A_{L,i}^{n_{L,i}}\tilde{\Delta}_{i}\subseteq \mu_{L,i}\tilde{\Delta}_{i}$~\cite{sasa2005}. 

Moreover, the following relationship can be inferred from \eqref{eq:estimation_error}
\begin{equation}
\begin{aligned}
x_i(k)&\in \hat{x}_i(k)\oplus  \tilde{\mathbb{S}}_{i}\,, \label{eq:estimation_error_set1}\\
\end{aligned}
\end{equation}
which yields a steady state assumption $\tilde{x}_i(0)\in\tilde{\mathbb{S}}_{i}$ because if the state estimation $\hat{x}_i$ lies in the tightened constraint set $\mathbb{X}_i\ominus \tilde{\mathbb{S}}_{i}$, the original state $x_i$ is guaranteed to lie in $\mathbb{X}_i$.

In this paper, the robust control policy is formed by an open-loop control solved by an optimization problem with nominal dynamics given by \eqref{eq:discrete_nominal_dynamics} and tightened constraints, and ancillary feedback control based on the observation of the state (from \eqref{eq:observer_dynamics}). The feedback controller is defined as:
\begin{equation}\label{eq:feedback_controller}
\begin{aligned}
u_i(k)=\bar{u}_i(k)+K_ie_i(k)\,,
\end{aligned}
\end{equation}
where $K_i$ is the gain of the feedback controller such that $A_{K,i}(=\! A_i\!+\!B_iK_i)$ meets the condition  $\lambda(A_{K,i})\!<\!1$. $e_i$ is the tracking error between the observer state and the nominal system, defined as:
\begin{equation}\label{eq:error_nominal}
e_i(k)=\hat{x}_i(k)-\bar{x}_i(k)\,.
\end{equation}
With the control action \eqref{eq:feedback_controller}, the closed-loop observer state satisfies:
\begin{equation}\label{eq:closeloop_observer}
\begin{aligned}
\hat{x}_i(k+1)=& A_i\hat{x}_i(k)+B_i\bar{u}_i(k)+B_iKe_i(k)\\&+L_{c,i}C_i\tilde{x}_i(k)+L_{c,i}\nu_i(k)\,,
\end{aligned}
\end{equation}
and the dynamics of the tracking error $e_i(k)$ are:
\begin{equation}\label{eq:error_nominal_dynamics}
\begin{aligned}
e_i(k+1)&=A_{K,i}e_i(k)+\bar{\delta}_i(k)\,,\\
\bar{\delta}_i(k)&=L_{c,i}C_i\tilde{x}_i(k)+L_{c,i}\nu_i(k)\,,\\
\bar{\delta}_i(k)&\in \bar{\Delta}_i= L_{c,i}C_i\tilde{\mathbb{S}}_i\oplus L_{c,i}\mathbb{V}_i\,.
\end{aligned}
\end{equation}
By analogy to the robust invariant tube for the observation error $\tilde{x}_i(k)$, the following conditions hold for the $e_i(k)$:
\begin{equation}\label{eq:error_nominal_2}
\begin{aligned}
e_i(k)&\in\bar{\mathbb{S}}_i, \forall k\,,\\
\hat{x}_i(k)&\in\bar{x}_i(k)\oplus\bar{\mathbb{S}}_i
\end{aligned}
\end{equation}
with the robust invariant tube for $e_i(k)$, 
$
\bar{\mathbb{S}}_i\!=\!(1\!-\!\mu_{K,i})^{-1}\oplus_{j=0}^{n_{K,i}-1}A_{K,i}^{j}\bar{\Delta}_i
$, where a finite integer $n_{K,i}\!\in\!\mathbb{N}_{>0}$ and a scalar $\mu_{K,i}\!\in\![0,1)$ satisfy $A_{K,i}^{n_{K,i}}\bar{\Delta}_{i}\subseteq \mu_{K,i}\bar{\Delta}_{i}$~\cite{sasa2005}, and $A_{K,i}\bar{\mathbb{S}}_i\oplus \bar{\Delta}_i\subseteq\bar{\mathbb{S}}_i$. 

Based on the definition of the estimation error \eqref{eq:estimation_error} and the tracking error \eqref{eq:error_nominal},
the actual state can be achieved by:
\begin{equation}\label{eq:actual_state}
\begin{aligned}
x_i(k)=\bar{x}(k)+e_i(k)+\tilde{x}_i(k)\,.
\end{aligned}
\end{equation}

The bound on the nominal state can be derived by combining the actual state and input constraints \eqref{eq:discrete_dynamics} and the feedback control policy \eqref{eq:feedback_controller},
\begin{subequations}\label{eq:nominal_bound}
\begin{align}
&\bar{x}_i(k)\in \bar{\mathbb{X}}_i\,,\label{eq:nominal_state_bound}\\
&\bar{u}_i(k)\in \bar{\mathbb{U}}_i\,,\\
&\bar{\mathbb{X}}_i= 
\left\{\begin{array}{lll} 
\mathbb{X}_i\ominus {\mathbb{S}}_i,\,\text{if }i=1 \,,\\
\mathbb{X}_i\ominus {\mathbb{S}}_i\ominus {\mathbb{S}}_h ,\,\text{if }i>1 \,,
\end{array}\right.
\nonumber \\
&\bar{\mathbb{U}}_i= \mathbb{U}_i\ominus (K_i \bar{\mathbb{S}}_i)\,, \nonumber
\end{align}
\end{subequations}
where ${\mathbb{S}}_i= \tilde{\mathbb{S}}_i\oplus\bar{\mathbb{S}}_i$, ${\mathbb{S}}_h= \tilde{\mathbb{S}}_h\oplus\bar{\mathbb{S}}_h$, and $\tilde{\mathbb{S}}_h$ and $\bar{\mathbb{S}}_h$ are the estimation error and tracking error derived invariant tubes of $\mathbb{X}_h$, respectively, which are generated to cope with the disturbance in the exchange information $\mathcal{I}(k_h)$ of the vehicle $h$ ($h\!\in\!\mathbb{N}(t_i(0))$, $h\!<\!i$, and $i\!\in\!\mathcal{C}_h\cup\mathcal{L}_h$). Moreover the following assumptions are imposed to guarantee the feasibility of \eqref{eq:nominal_bound}:
\[
\begin{aligned}
{\mathbb{S}}_i\subset \mathbb{X}_i, \quad
K_i\bar{\mathbb{S}}_i\subset \mathbb{U}_i\,.
\end{aligned}
\]
The assumptions above ensure the existence of tightened sets for the nominal state $\bar{x}_i(k)$ and input $\bar{u}_i(k)$
such that the actual state $x_i(k)$ and input $u_i(k)$ of the controlled system $x_{i}(k+1)\!=\! A_ix_{i}(k)\!+\!B_iu_{i}(k)\!+\!\omega_{i}(k)$ with feedback-loop $u_i(k)\!=\!\bar{u}_i(k)\!+\!K_i(\hat{x}_i(k)-\bar{x}_i(k))$ satisfy the original constraints $x_i(k)\in\mathbb{X}_i$ and $u_i(k)\in\mathbb{U}_i$ for all admissible disturbances $\omega_i(k)$ and $\nu_i(k)$.

Therefore, this work proposes an RMPC-based HRCS designed based on the nominal system \eqref{eq:nominal_dynamics}.
At an update instant $k$, the tube-based RMPC in the lower-level finds the optimal control sequence $\bar{\mathbf{u}}^*_{i}(k)\!=\!\{\bar{u}^*_{i}(k|k),\bar{u}^*_{i}(k+1|k),\ldots,\bar{u}^*_{i}(k\!+\!N_p\!-\!1|k)\}$ and the optimal state sequence $\bar{\mathbf{x}}^*_{i}(k)\!=\!\{\bar{x}_i^*(k|k),\bar{x}^*_{i}(k+1|k),\ldots,\bar{x}^*_{i}(k+N_p|k)\}$ by solving the following convex RMPC problem:

\begin{problem}\label{prob:ocp4}
\begin{subequations}\label{eq:tubempc}
\begin{align}
&\min_{\bar{\mathbf{u}}_{i}(k)}\hspace{2mm} \bar{J}_{d,i}(\bar{\mathbf{x}}_{i}(k),\bar{\mathbf{u}}_{i}(k))\\
&\textbf{s.t.}\hspace{2mm}   \bar{x}_{i}(k|k)=\hat{x}(k)\\
&\hspace{6.5mm} \bar{x}_{i}(k+j+1|k) \!=\! A_i\bar{x}_{i}(k+j|k) \!+\!B_i\bar{u}_{i}(k+j|k)\\
&\hspace{6.5mm} \bar{u}_i(k+j|k)\in\bar{\mathbb{U}}_i\\
&\hspace{6.5mm} \bar{x}_i(k+j|k)\in\bar{\mathbb{X}}_i,\quad j=0,1,\cdots,N_p-1\\
&\textbf{given}:\, \hat{x}_{i}({0})\in
\left\{\begin{array}{lll} 
 \mathbb{X}_i\ominus\mathbb{S}_i,\,\text{if }i=1 \\
 {\mathbb{X}_i\ominus\mathbb{S}_i\ominus\mathbb{S}_h ,\,\text{if }i>1 }
\end{array}\right.
\end{align}
\end{subequations}
\end{problem}
where $\bar{\mathbf{{x}}}_{i}(k)\!\in\!\mathbb{R}^{(N_p+1)\times 2}$, $\bar{\mathbf{{u}}}_{i}(k)\!\in\!\mathbb{R}^{N_p\times 2}$, $k{=} 0,1,\ldots,\alpha$, $\hat{x}(k)$ is obtained by the observer \eqref{eq:closeloop_observer} at the instant $k$, and the given and bounded $\hat{x}(0)$ is the measured initial state (entry state at CZ) satisfying Assumptions \ref{ass:initial}. Moreover, $\bar{J}_{d,i}(\bar{\mathbf{x}}_{i}(k),\bar{\mathbf{u}}_{i}(k))$ is the augmented objective function of \eqref{eq:J_decentralized}, which will be designed next.

Due to the receding horizon nature of MPC, there is a potential infeasibility issue when vehicles reach the MZ from perpendicular or opposite directions, and the constraints related to MZ are not initially considered.
To prevent the infeasibility, the objective function \eqref{eq:J_decentralized} in the lower-level RMPC of HRCS for vehicle $i\!\in\!\{\mathcal{L}_h,\mathcal{O}_h\}$ is augmented by relaxing constraints \eqref{eq:lateral} and \eqref{eq:opposite} as follows:
\begin{multline}\label{eq:decentralized_objective}
\bar{J}_{d,i} = {J}_{d,i}
+W_3\,(\max\{0,\Delta t_i^{\mathcal{L}}-\Gamma_i^{\mathcal{L}}\})^2
\\+W_4\,(\max\{0,\Delta t_i^{\mathcal{O}}-\Gamma_i^{\mathcal{O}}\})^2\,,
\end{multline}
where $\Delta t_i^{\mathcal{L}},\,\Delta t_i^{\mathcal{O}}\in\mathbb{R}_{>0}$ are predefined tuneable parameters, and 
$\Gamma_i^{\mathcal{L}}$ and $\Gamma_i^{\mathcal{O}}$ are the time difference of the constraints \eqref{eq:lateral} and \eqref{eq:opposite},
\begin{equation}\label{eq:si_lateral}
\Gamma_i^{\mathcal{L}}=t_{i}^\mathcal{L} - \hat{t}_{h}^{\mathcal{L}}\,,
\end{equation}
\begin{equation}\label{eq:si_opposite}
\Gamma_i^{\mathcal{O}}=t_{i}^{\mathcal{O}} - \hat{t}_{h}^{\mathcal{O}}\,.
\end{equation}
Tuneable parameters $\Delta t_i^{\mathcal{L}}$ and $\Delta t_i^{\mathcal{O}}$ of a similar nature have been introduced in~\cite{riegger2016centralized}, but in a rather simpler non-robust framework for autonomous intersection management with a limited number of vehicles, simple vehicle modeling, and a different shrinking MPC horizon technique.
Thus, if the time difference between vehicle $i$ and $h$ is large enough ($\Gamma_i^{\mathcal{L}}>\Delta t_i^{\mathcal{L}}$ and $\Gamma_i^{\mathcal{O}}>\Delta t_i^{\mathcal{O}}$), for $i\in \mathcal{L}_h\cup \mathcal{O}_h$ the augmented objective function \eqref{eq:decentralized_objective} is identical to the original objective function \eqref{eq:J_decentralized}. On the other hand, the augmented terms in \eqref{eq:decentralized_objective} encourage the $i$th vehicle to gradually accelerate or decelerate in advance to maintain a time gap equal to the predefined parameters from the $h$th vehicle over the entire mission ($\Gamma_i^{\mathcal{L}}=\Delta t_i^{\mathcal{L}}$ and $\Gamma_i^{\mathcal{O}}=\Delta t_i^{\mathcal{O}}$). In this way, the hard braking to satisfy the constraints \eqref{eq:lateral} and \eqref{eq:opposite} can be avoided when the $i$th vehicle enters the MZ.
Taking $\Gamma_i^{\mathcal{L}}$ as an example, it is specified as below:
\begin{subequations}\label{eq:decentralize_lateral_right}
\begin{align}
& \text{1) if } k+N_p<\alpha_1\wedge k_h-k\leq \alpha_2+\alpha_3 \nonumber \\
& t_i^{\mathcal{L}} = t_i(k+N_p|k) \nonumber \\
&\hat{t}_{h}^{\mathcal{L}}\!=\! 
\left\{\begin{array}{lll} 
\hspace{-3mm} \displaystyle\,\,\bar{t}_{h}\!(k_h\!+\!l\!+\!1|k_h),\,\displaystyle\text{if}\,\, l\!<\!N_p,\,k_h\!+\!l\!+\!1\!=\!k\!+\!\alpha_2\!+\!\alpha_3\!, \vspace{2mm}\\
\hspace{-3mm} \vspace{1mm}\displaystyle \,\,\bar{t}_{h}(k_h\!+\!N_p|k_h)\!+\!\frac{[k\!+\!\alpha_2\!+\!\alpha_3-\!(k_h\!+\!N_p)]\Delta s}{\bar{v}_{h}(k_h+N_p|k_h)},\\
\hspace{30mm}\displaystyle \text{if}\,\,{k}_h\!+\!N_p<k\!+\!\alpha_2\!+\!\alpha_3,
\end{array}\right.  \label{eq:decentralize_lateral_right1}\\
& \text{2) if } k+N_p\geq\alpha_1 \wedge k<\alpha_1  \nonumber \\
& t_i^{\mathcal{L}} = t_i(k+j+1|k),\,k+j+1=\alpha_1 \nonumber \\
&\hat{t}_{h}^{\mathcal{L}}\!=\! 
\left\{\begin{array}{lll} 
\hspace{-3mm} \displaystyle\,\,\bar{t}_{h}\!(k_h\!+\!l\!+\!1|k_h),\,\displaystyle\text{if}\,\, l\!<\!N_p,\,k_h\!+\!l\!+\!1\!=\!\alpha\!, \vspace{2mm}\\
{ \hspace{-3mm} \vspace{1mm}\displaystyle \,\,\bar{t}_{h}(k_h\!+\!N_p|k_h)\!+\!\hspace{-4mm}
\sum_{\beta=1}^{\alpha\!-\!(k_h\!+\!N_p)}\hspace{-2mm}\frac{\Delta s}{\hat{v}(\beta)},}\\
\hspace{30mm}\displaystyle \text{if}\,\,{k}_h\!+\!N_p<\alpha,
\end{array}\right.
\label{eq:decentralize_lateral_right2}
\end{align}
\end{subequations}
where $l\!\in\!\mathbb{N}_{[0,N_p-1]}$, and $k_h$, similarly to the earlier definition, stands for the associated step of the ${h}$th CAV such that $t_{h}(k_h \Delta s)\!=\!t_{i}(k \Delta s)$.
The first case activates when the $i$th and the $h$th vehicles have a similar remaining distance to the entry of the MZ. The formulation encourages the $i$th vehicle to preserve a predefined time gap $\Delta t_i^{\mathcal{L}}$ for $\alpha_2\!+\!\alpha_3$ step-ahead space (equivalent to distance length $S\!+\!l_h$) to the $h$th vehicle before entering the MZ. As such, the effect of the lateral collision constraint \eqref{eq:lateral} is considered from the start of the mission, and eventually the potential infeasibility issue in MPC framework when the $i$th vehicle and the $h$th vehicle reach the MZ at a similar time can be avoided.
In the second case, as reflected by \eqref{eq:decentralize_lateral_right2}, the lateral collision avoidance constraint~\eqref{eq:lateral} should still remain, and in discretized format it is described as follows:
\begin{equation}
    \label{eq:lateral_decentralize_format}
{t}_{i}^{\mathcal{L}}-\hat{t}_{h}^{\mathcal{L}}\geq0, \quad\text{if }\, k+N_p\geq\alpha_1 \wedge k<\alpha_1,
\end{equation}
where the upper part of $\hat{t}_h^{\mathcal{L}}$ in \eqref{eq:decentralize_lateral_right2} implies that the time of the $h$th vehicle at the position ($L_h$) can be obtained straightforwardly in the present space window. If the desired time at the position ($L_h$) is unavailable, the lower part of $\hat{t}_h^{\mathcal{L}}$ in \eqref{eq:decentralize_lateral_right2} provides a conservative estimation method, { where the ${h}$th CAV will finish the rest of the mission at a conservatively (the slowest) estimated speed $\hat{v}(\beta)$, which is calculated based on the terminal speed of the present horizon and the minimum deceleration of the kinetic energy as below:
\begin{equation*}
\begin{aligned}
&\hat{v}(\beta)\!=\! \\
&\max\left\{{\sqrt{2\left(\bar{E_h}(k_h\!+\!N_p|k_h)\!-\beta \left.\frac{dE_h(s)}{ds}\right|_{\min}\right)/m_i}},v_{\min}\right\} 
\end{aligned}
\end{equation*}
with
\begin{equation*}
\begin{aligned}
\left.\frac{dE_h(s)}{ds}\right|_{\min} = F_{w,\min}-F_{r,i}-f_{d,i}v_{\max}^2\,\,\,\text{being a constant.}
\end{aligned}
\end{equation*}
Note that the estimated speed $\hat{v}(\beta)$ and the index $\beta\in[1,\alpha\!-\!(k_h\!+\!N_p)]$ can be predetermined before the optimization.
}
By analogy to the definition \eqref{eq:decentralize_lateral_right} of $\Gamma_i^{\mathcal{L}}$, the definition and estimation of $\Gamma_i^{\mathcal{O}}$ can also be obtained. 
The prioritization constraint \eqref{eq:opposite} in discretized form is then written as follows:
\begin{equation}
    \label{eq:opposite_decentralize_format}
{t}_{i}^{\mathcal{O}}-\hat{t}_{h}^{\mathcal{O}}\geq0, \quad\text{if }\, k+N_p\geq\alpha_1 \wedge k<\alpha_1.
\end{equation}

{
\begin{corollary}
By finely tuning the parameters  $\Delta t_i^{\mathcal{L}}$ and $\Delta t_i^{\mathcal{O}}$, the solution of the MPC problem always satisfies condition \eqref{eq:syst3} with equality, and therefore the solution of the MPC problem is valid under \cref{prop:tightness}.
\end{corollary}

\textit{Proof:} Assume a feasible solution sequence is found within an MPC horizon with $k+j+1=\alpha_1$ such that:
\begin{subequations}
\begin{align}
&\bar{E}_i^*(k+j+1|k)=
\bar{E}_i^*(k|k)+\sum_{j_k=0}^{j}\{\bar{F}^*_{w,i}(k+j_k|k) \\
&\hspace{3.5cm}-F_{r,i}+\varepsilon_i \bar{E}_i^*(k+j_k|k)\}\Delta s \nonumber\\
&\bar{t}_i^*(k+j+1|k)=\bar{t}_i^*(k|k)+\sum_{j_k=0}^{j}\{\bar{\zeta}^*_{i}(k+j_k|k) \}\Delta s\\
&\bar{\zeta}^*_{i}(k+j_k|k) \geq \frac{1}{\sqrt{2\bar{E}_i^*(k+j_k|k)/m_i}} \label{eq:corollary_tightness}\\
& \bar{t}_{i}^*(k+j+1|k) - \hat{t}_{h}^{\mathcal{L}} \geq 0, \label{eq:corollary_lateral}\\
& \bar{t}_{i}^*(k+j+1|k) - \hat{t}_{h}^{\mathcal{O}} \geq 0, 
\\ &\textbf{given: }\bar{E}_i^*(k|k)= {E}_i^*(k),\,\bar{t}_i^*(k|k)={t}_i^*(k).
\end{align}
\end{subequations}

Considering $\bar{\zeta}_i^*\!=\!\bar{\breve{\zeta}}_i+\Delta \zeta_i$, where $\bar{\breve{\zeta}}_i$ is the solution in \eqref{eq:breve_zeta} that satisfies \cref{prop:tightness} with $\Delta \zeta_i>0$, the constraint \eqref{eq:corollary_lateral} can be expanded as:
\vspace{-.1cm}
\begin{align}
0&\leq  \bar{t}_{i}^*(k+j+1|k) - \hat{t}_{h}^{\mathcal{L}}  \nonumber\\
&\leq \bar{t}_i^*(k)+\sum_{j_k=0}^{j}\{\bar{\zeta}_{i}^*(k+j_k|k)\}\Delta s-\hat{t}_h^{\mathcal{L}} \nonumber\\
&\leq \bar{t}_i^*(k)+\sum_{j_k=0}^{j}\{\bar{\breve{\zeta}}_{i}(k+j_k|k)+\Delta \zeta_i(k+j_k|k) \}\Delta s-\hat{t}_h^{\mathcal{L}} \nonumber\\ 
& \leq \left[\bar{t}_i^*(k)+\sum_{j_k=0}^{j}\{\Delta \zeta_i(k+j_k|k)\}\Delta s \right] \nonumber\\
&\hspace{3.3cm} +\sum_{j_k=0}^{j} \bar{\breve{\zeta}}_{i}(k+j_k|k)\Delta s-\hat{t}_h^{\mathcal{L}}\,.
\label{eq:corollary_lateral_expand}
\end{align}

In virtue of the augmented objective function \eqref{eq:decentralized_objective}, it can be claimed that if there exists a feasible solution, by finely tuning $\Delta t_i^{\mathcal{L}}$ to enforce the CAV $i$ to accelerate at a slow rate or start to decelerate before step $k$, the feasible solution with initial condition $\bar{\breve{t}}_{i}(k)$ at step $k$ can be found such that: 
\begin{equation}\label{eq:corollary_lateral_expand1}
 \bar{\breve{t}}_{i}(k)\!\geq\! \bar{t}_i^*(k)\!+\!\sum_{j_k=0}^{j}\{\Delta \zeta_i(k+j_k|k)\}\Delta s\,,
\end{equation}
and constraint \eqref{eq:corollary_lateral_expand} can be reformulated as
\begin{equation}\label{eq:corollary_lateral_expand2}
\bar{\breve{t}}_{i}(k) + \sum_{j_k=0}^{j} \bar{\breve{\zeta}}_{i}(k+j_k|k)\Delta s-\hat{t}_h^{\mathcal{L}} \geq 0\,.
\end{equation}
Thus, the feasible solution with its initial condition satisfying \eqref{eq:corollary_lateral_expand1} and \eqref{eq:corollary_lateral_expand2} consequently satisfies the tightness of \eqref{eq:syst3}.
By analogy to the forgoing analysis of \eqref{eq:corollary_lateral}, the same design approach can be applied to $\Delta t_i^{\mathcal{O}}$.
As a result, by suitably selecting the value of parameter $\Delta t_i^{\mathcal{L}}$ and $\Delta t_i^{\mathcal{O}}$, the solution of the decentralized MPC framework for autonomous intersection can be guaranteed to be valid under \cref{prop:tightness}.
}

\begin{remark}
The convex tightened sets $\bar{\mathbb{U}}_i$ and $\bar{\mathbb{X}}_i$ in \eqref{eq:nominal_bound} for each vehicle $i$ can be computed offline in a decentralized manner to increase computational efficiency. 
\end{remark}

\begin{remark}
Given a crossing order obtained in the upper-level of the HRCS, there might be a case where no valid solution can be found in the lower-level due to the discrepancies appearing between the upper and lower optimization problems. This can be addressed by recursively solving the lower-level problem with continually reducing $v_{\min}$, and it terminates when a valid { feasible} solution is found. 
\end{remark}

\section{Numerical Results}
\label{sec:simulation}
{ The evaluation of the proposed HRCS is listed below: 
1) the effectiveness of the proposed convexified HRCS is verified and compared with benchmark solutions obtained by a nominal MPC-based strategy under the same initial conditions;
2) a comparison between the HRCS and a benchmark using the same tube-based MPC algorithm but following the FIFO policy to show the better performance of the designed upper-level crossing order scheduler;
3) the lower-level tube-based RMPC \eqref{eq:tubempc} is solved for different weighting combinations $\{W_1, W_2\}$ of \eqref{eq:decentralized_objective} under fixed values of $\{W_3,W_4\}$ and a series of different arrival rates to show the trade-off between energy cost and travel time as well as the impact of the traffic density on the overall optimality;
4) the trade-offs between the robustness and the optimality
are investigated;
5) the computational time with respect to the size of the sampling interval is investigated to show the validity of the method in potential practical application.
}

{%
For the sake of fair comparison, a terminal speed condition is imposed for all CAVs: 
\begin{equation}
\frac{1}{2}mv_f^2 -\gamma_i \leq E_i(L_i) \leq \frac{1}{2}mv_f^2 +\gamma_i \label{eq:terminalv}
\end{equation}
where $v_f \!\in\![v_{\min},v_{\max}]$ is a predefined terminal speed, and $\gamma_i\!\in\!\mathbb{R}_{>0}$ is an auxiliary optimization variable,  which is introduced to avoid potential feasibility issues caused by additive disturbances (i.e. process noise $\omega_{i}$ and measurement noise $\nu_{i}$). The terminal condition \eqref{eq:terminalv} encourages all CAVs to reach the same terminal speed when leaving the MZ by minimizing $\gamma_i$ in the objective function, $\bar{J}_{d,i} + W_5\cdot \gamma_i$. Note that it is straightforward to relax \eqref{eq:terminalv} with non-uniform terminal speeds. 
}

For the sake of simplicity, the simulation in this work assumes all CAVs to be identical, with
the main characteristic parameters of each vehicle model summarized in Table~\ref{tab:vehicledata}. 
In the following case studies, the parameters of the intersection are $L \!=\! 150$~m and $S \!=\! 10$~m with sampling interval $\Delta s\!=\!2$~m. The bounds of the external disturbances $\omega_i$ and $\nu_i$ are
$
\overline{\omega}_{E,i}\!=\!
\frac{1}{2}m\overline{\omega}_{v,i}^2, \,
\overline{\nu}_i\!=\![
\frac{1}{2}m\overline{\nu}_{v,i}^2,\,
\overline{\nu}_{t,i}
]^{\top}
$
where $\overline{\omega}_{v,i}\!=\!\overline{\nu}_{v,i}\!=\!0.1~m/s$ and $\overline{\nu}_{t,i}\!=\!0.1~s$.
The values of the tuneable parameters in \eqref{eq:decentralized_objective} are predefined to be $\Delta t_i^{\mathcal{L}}\!=\!\Delta t_i^{\mathcal{O}}\!=\!0.4$~s based on consideration of the intersection dimensions.
The minimum allowed time gap in \eqref{eq:TTCapprox} is set to $t_{\delta}\!=\!1$~s, which is chosen as such due to the expectation that CAVs can cope safely with further reduced headway distances than human-driven vehicles.
For illustrative and comparative purposes, let us consider a case with $N\!=\!20$ CAVs all of which are assumed to leave the intersection (the MZ) at the same terminal speed $v_f\!=\!10$~m/s. Without loss of generality, the control problem is initialized with randomized initial conditions $v_i(0)$ and $t_i(0)$ for all CAVs subject to the constraints imposed in 
Assumptions~\ref{ass:initial}. In particular, CAVs’ initial speeds follow a uniform distribution within $[v_{\min},\, v_{\max} ]$, while their arrival times, $t_i(0)$, follow a Poisson distribution. Moreover, the entry direction of each CAV is also randomly generated. 
The proposed HRCS is solved by using YALMIP and MOSEK~\cite{Lofberg2004} in Matlab on a personal computer with Intel Core i5 2.9 GHz and 8 GB of RAM.

\begin{figure}[t!]
\centering
\vspace{-3mm}
\includegraphics[width=.95\columnwidth]{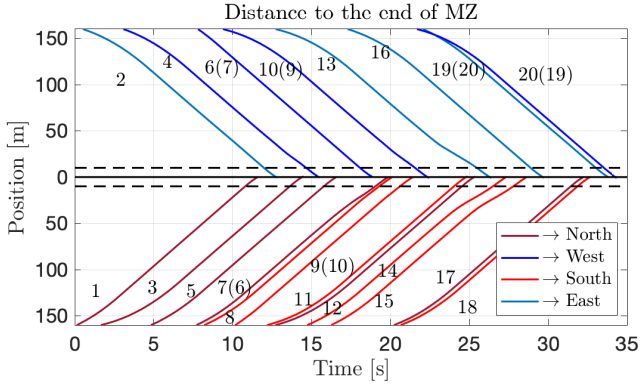}
\\[-2ex]
\caption{{ Traveled distance trajectories (distance to the end of MZ) by solving the HRCS  subject to an average travel time $12.31$~s at an arrival rate of 800 veh/h per lane and with prediction horizon length $N_p\!=\!15$. The horizontal dashed lines correspond to the entry of the MZ, while the horizontal continuous black line denotes the end of the MZ. The four vehicle heading directions are denoted using different colors. Note that the numbers in the brackets highlight the arriving orders of the vehicles at the CZ, which are different from their order entering the MZ. The upper-level scheduler sorts the vehicles in order of $\mathcal{N}=\{1,2,3,4,5,7,6,8,10,9,11,12,13,14,15,16,17,18,20,19\}$. }}\vspace{-1mm}
\label{fig:trajectories}
\end{figure}

{
In the first instance, the proposed method is solved with prediction horizon length $N_p\!=\!15$ at an arrival rate of 800 veh/h (vehicles per hour) per lane, which is ordinary for practical intersections. The weighting factors are set to emphasize more on the travel time rather than energy consumption in the objective function. Fig.~\ref{fig:trajectories} presents the optimal traveled distance profile subject to an average travel time $12.31$~s. As it can be seen, by giving the order $\mathcal{N}$ defined in the upper-level, the lower-level RMPC coordinates the CAVs such that no CAVs violate the rear-end and lateral collision constraints, which verifies the validity of the optimal solution. An example of lateral collision avoidance can be found in vehicle 12 and vehicle 13, where vehicle 13 will not be allowed to enter the MZ until the ahead vehicle 12 has left. 
The effectiveness of rear-end collision avoidance can be identified as the solution has no intersections between trajectories of the same color throughout the CZ. In addition, the average, minimum and maximum traveling times of all CAVs in Fig.~\ref{fig:trajectories} are also provided in Table.~\ref{tab:statistic}.

\begin{table}[!t]
\centering 
\caption{Electric vehicle model parameters}
\label{tab:vehicledata} 
\begin{tabular*}{1\columnwidth}{l @{\extracolsep{\fill}} cl}
\hline
\hline
 symbol & value & description\\
 \hline
 $m_i$ & 1200\,kg & vehicle mass \\
 $l_i$ & 4\,m & vehicle body length\\
 $r_{w,i}$ & 0.3\,m & wheel radius\\
 $g_{r,i}$ & 3.5 & transmission gear ratio\\
 $f_{r,i}$ & $0.01\,$ & rolling resistance coefficient\\
 $f_{d,i}$ & $0.47\,$ & air drag resistance coefficient\\
 $v_{\min}$    &$0.1$\,m/s& minimum velocity \\
 $v_{\max}$    &$15$\,m/s& maximum velocity \\
 $a_{\min,i}$    & -6.5\,m$\cdot$s$^{-2}$ & minimum acceleration \\
\hline
\hline
\end{tabular*}
\end{table}

}

\begin{figure}[t!]
\centering
\includegraphics[width=1\columnwidth]{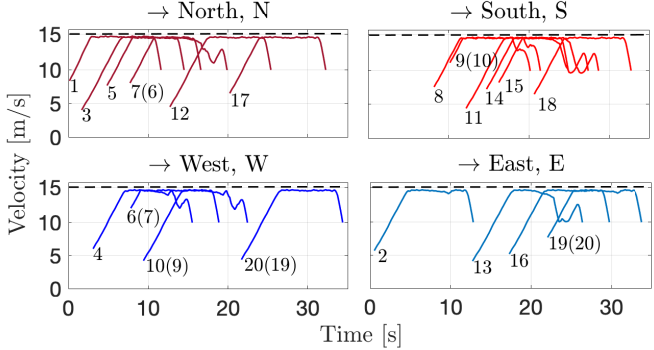}
\\
[-2ex]
\caption{Optimal speed profiles by solving the HRCS subject to an average travel time $12.31$~s for all CAVs at an arrival rate of 800 veh/h per lane. Note that the numbers in the brackets highlight the arriving orders of the vehicles at the CZ, which are different from their order entering the MZ.}
\label{fig:RMPC_velocityprofile}
\end{figure}

\begin{figure}[t!]
\centering
\includegraphics[width=1\columnwidth]{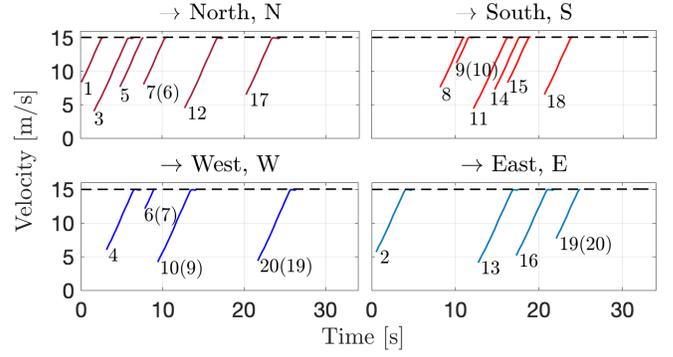}
\\
[-2ex]
\caption{Optimal speed profiles by solving the nominal MPC-based strategy under the same initial conditions as the results in Fig.~\ref{fig:RMPC_velocityprofile}. Note that the numbers in the brackets highlight the arriving orders of the vehicles at the CZ, which are different from their order entering the MZ.}
\label{fig:Nominal_velocityprofile}
\end{figure}

{ To verify the robustness of the proposed method, a comparison of the optimal speed trajectories with benchmark solutions obtained by a nominal MPC-based strategy under the same initial conditions and an average travel time $12.31$~s can be found in Fig.~\ref{fig:RMPC_velocityprofile} and Fig.~\ref{fig:Nominal_velocityprofile}, respectively, where the speed profiles in both scenarios are grouped based on the heading directions for illustration purposes.
}
As it can be seen in Fig.~\ref{fig:RMPC_velocityprofile}, the velocity profiles of most CAVs solved by the proposed robust method tend to cruise at a constant speed and apply more intensive braking when approaching the exit of the MZ, which is due to the energy recovery in the powertrain of the battery electric vehicles. In some cases, if two vehicles have a potential collision inside the MZ and reach the MZ at a close time, the speed may not follow the foregoing trajectories due to the introduction of the augmented objective function \eqref{eq:decentralized_objective}. For example, the speed profiles of CAV 7 and CAV 13 exhibit early deceleration when approaching the MZ to preserve enough time gaps for lateral collision avoidance.
To deal with the impact of external disturbances, 
the robust invariant tubes in HRCS prevent the peak speed from reaching the maximum allowed velocity by a margin even though in this case the emphasis is more on travel time than energy consumption.
However, the peak speed in the nominal MPC-based strategy case (see in Fig.~\ref{fig:Nominal_velocityprofile}) remains at the constant value at $v_{\max}=15$~m/s, which leads to fragile feasibility when it comes to disturbances. As a result, infeasible solutions are yielded by the nominal MPC-based strategy, which violates the speed constraints.

\begin{table}[t!]
\centering 
\caption{{ Travel time and velocity performances of all CAVs by solving the HRCS at an arrival rate of 800 Veh/h per lane subject to average travel time at 12.31 s and with $N_p\!=\!15$ }.}
\label{tab:statistic} 
\begin{tabular*}{1\columnwidth}{c @{\extracolsep{\fill}} c@{\extracolsep{\fill}} c@{\extracolsep{\fill}}c@{\extracolsep{\fill}}c}
\hline
\hline
  & Min & Max & Average & Variance\\
 \hline
 Travel Time [s] & 11.04 & 13.65 & 12.26 & 0.3934 \\
 Velocity [m/s]  & 3.89 & 14.82 & 13.53 & 4.6022\\
\hline
\hline
\end{tabular*}
\end{table}


%

%

\begin{figure}[t!]
\centering
\includegraphics[width=.95\columnwidth]{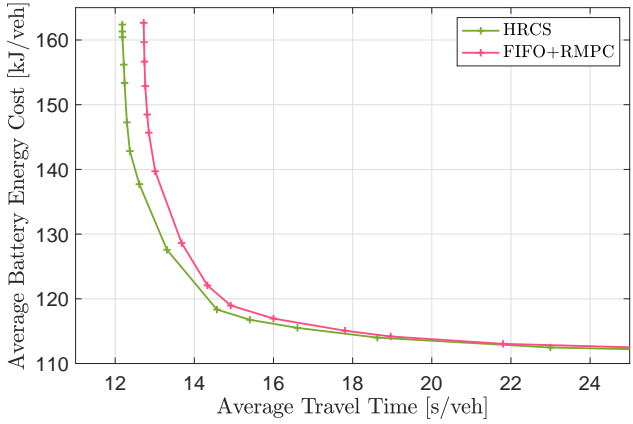}\\[-2ex]
\caption{{ Comparison of the energy-time cost trade-off between the proposed method and a benchmark using the same tube-based MPC following the FIFO policy with prediction horizon length $N_p\!=\!15$ at an arrival rate of 800 veh/h per lane.}}
\label{fig:trade_off_FIFO}
\end{figure}

{
The optimality of the proposed method is firstly investigated by comparing with a benchmark solution using the same decentralized tube-based MPC algorithm for the lower-level trajectory optimization but following the FIFO policy. As it can be seen in Fig.~\ref{fig:trade_off_FIFO}, owing to the proposed upper-level crossing scheduler, the HRCS as compared to the benchmark can reach a more energy-time efficient solution. In particular, when the average travel time is 12.72 s, the proposed method can save up to 16.35\% energy consumption over the benchmark. 
}

\begin{figure}[htb!]
\centering
\vspace{-3mm}
\includegraphics[width=.95\columnwidth]{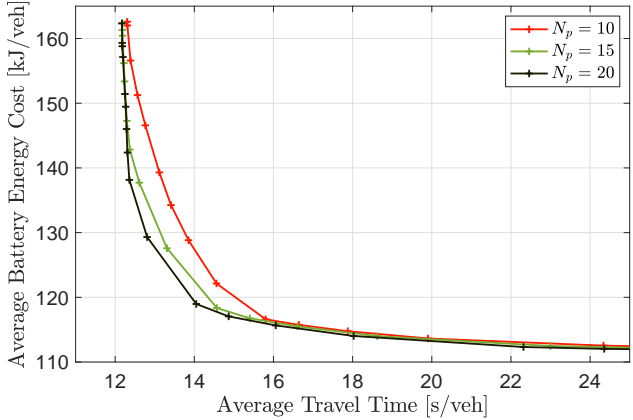}\\[-2ex]
\caption{The trade-offs between average battery energy consumption and average travel time at the arrival rate of 800 veh/h per lane for the decentralized HRCS with varied prediction horizon length $N_p= \{10,15,20\}$.}
\label{fig:trade_off_horizon}
\end{figure}

\begin{figure}[t!]
\centering
\vspace{-3mm}
\includegraphics[width=.95\columnwidth]{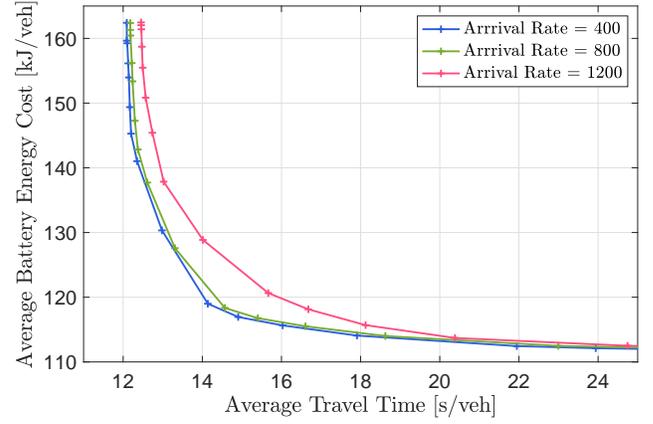}\\[-2ex]
\caption{The trade-offs between average battery energy consumption and average travel time with prediction horizon length $N_p\!=\!15$ at different arrival rates from 400 to 1200 veh/h per lane for the decentralized HRCS.}
\label{fig:trade_off_arrivalrate}
\end{figure}

In order to investigate the impact of the prediction horizon length and the traffic density in HRCS, the trade-off between travel time and energy consumption for a series of combinations of the weight factors, $W_1$ and $W_2$, (under the same initial conditions and fixed $W_3$, $W_4$ and $W_5$) and for different prediction horizon lengths and arrival rates, are presented in Fig.~\ref{fig:trade_off_horizon} and Fig.~\ref{fig:trade_off_arrivalrate}, respectively. For both cases, the Pareto front results point out the importance of examining the energy-time trade-off, as a small change in the travel time can significantly affect the energy efficiency. For example, when the proposed HRCS method is used with prediction horizon length $N_p= 15$ under the arrival rate of 800 veh/h per lane, an increase of 20\% in travel time (from 12.18~s to 14.6~s) can lead to an average energy consumption reduction of 26.71\%, while further increase in travel time can eventually yield up to 30.43\% energy consumption reduction. In Fig.~\ref{fig:trade_off_horizon}, the comparison among the three prediction horizon lengths indicates that the overall optimality increases as the $N_p$ increases. The reason is that increasing the horizon tends to enhance the ability to anticipate the future behavior of each CAV to satisfy the collision constraints, at a price of a higher computational burden. However, the improvement in terms of energy consumption optimality from $N_p\!=\! 15$ to $20$ is less than 6.1\% for most cases, and therefore, the subsequent studies in this paper adopt $N_p\!=\!15$ in the lower-level RMPC of the proposed HRCS.

Moreover, by comparing the results at different arrival rates from 400 to 1200 veh/h per lane, it can be observed in~Fig.\ref{fig:trade_off_arrivalrate} that the overall optimality deteriorates as the traffic density increases. The reason is that a higher arrival rate implies a higher traffic density condition, where the motions of vehicles are more restrained by the surrounding vehicles, and therefore, the optimal solution tends to be compromised by collision avoidance requirements. When most weight is placed on the travel time term ($W_2\!\gg\! W_1$), an optimality gap can be observed between the cases with the arrival rate of 800 and 1200 veh/h per lane, which indicates that with an emphasis on the travel time minimization, the optimization encourages the CAVs to travel at maximum speed, which yields more restrictive solutions due to the tougher collision avoidance constraints in such cases, and the restrictiveness rises as the arrival rate increases.
As the weight $W_1$ for the energy cost is gradually increased, the optimality deteriorates, resulting in a maximum gap of 12.4\% for the case with the arrival rate of 800 veh/h per lane as compared to the case with the arrival rate of 1200 veh/h per lane when the average traveled time is 12.45~s, and after this time the gap gets closer and gradually becomes negligible. The reason is that when the travel time is relaxed, there exists more room for speed optimization, and optimal solutions in terms of energy consumption become similar.
Finally, it has been found that a further decrease in the arrival rate below 400 veh/h makes a negligible impact on the optimality, as the traffic is sufficiently sparse to allow free optimization of each velocity trajectory without concession to other vehicles.

To construct the relationship between the average time gap and energy consumption, an example case of the solutions caused by tougher collision constraints (with more emphasis on the traveled time) is investigated. Table.~\ref{tab:KPI} presents the average time gap at different arrival rates with a fixed average traveled time 12.51 s. As it can be seen, there is an upward trend from left to right (as the arrival rate increases from 400 veh/h to 1200 veh/h) in the energy cost from 139.10 kJ to 164.48 kJ. Meanwhile, the average time gap decreases from 9.29 s to 1.69 s as the arrival rate increases from 400 veh/h to 1200 veh/h. This can be understood that increased traffic density could result in severe congestion and more acceleration/deceleration behavior, and therefore reduced time gap and higher energy consumption.

\begin{table}[t!]
\centering 
\caption{Average time gap at different arrival rates with a fixed average travel time 12.51 s with $t_{\delta} \!=\! 1$~s.}
\label{tab:KPI} 
\begin{tabular*}{1\columnwidth}{l @{\extracolsep{\fill}} c@{\extracolsep{\fill}} c@{\extracolsep{\fill}}c}
\hline
\hline
 Arrival Rate [veh/h] & 400 & 800 & 1200\\
 \hline
 Average Energy Cost [kJ] & 139.10 & 139.40 & 164.48 \\
 Average Time Gap [s] & 5.44 & 2.57 & 1.69\\
 Minimum Time Gap [s] & 3.30 & 1.26 & 1.10 \\
 Maximum Time Gap [s] & 9.29 & 4.56 & 2.98\\
\hline
\hline
\end{tabular*}
\end{table}

{
To further investigate the trade-off between robustness and optimality, further simulation trials are carried out, where disturbance bounds are conservatively used in the design of RMPC due to the lack of precise knowledge of the disturbances (commonly encountered in the practice). Recalling the disturbance bounds given in Section~\ref{sec:simulation} (identical for all disturbance sources $\overline{\omega}_{v,i}\!=\!\overline{\nu}_{v,i}\!=\!0.1~m/s$ and $\overline{\nu}_{t,i}\!=\!0.1~s$), under the same initial conditions and disturbances, two additional RMPCs are designed and simulated with more conservative bounds of 0.2 and 0.3, respectively. The comparative results are shown in Table.~\ref{tab:KPI_noise}. 
As it can be seen, when the average energy cost is fixed at 146.55~kJ, doubling the bound in the RPMC design can increase the travel time by 4.72\%, and the figure goes up to 10.8\% when the bound is tripled. The results can be understood that the greater disturbance bound leads to the more conservative RMPC design in terms of the constraints tightening (reduced feasibility), which is reflected in the minimum time gaps shown in the Table. Nevertheless, the maximum time gaps are not affected as no upper limit is imposed for the time gap between CAVs.
}
\begin{table}[t!]
\centering 
\caption{Average travel time and time gap with a fixed average energy cost 146.55 kJ and $t_{\delta} \!=\! 1$~s for different RMPC designs.}
\label{tab:KPI_noise} 
\begin{tabular*}{1\columnwidth}{l @{\extracolsep{\fill}} c@{\extracolsep{\fill}} c@{\extracolsep{\fill}}c}
\hline
\hline
 $\{\bar{\omega}_{v,i},\bar{\nu}_{v,i},\bar{\nu}_{t,i}\}$ & \{0.1,\,0.1,\,0.1\} & \{0.2,\,0.2,\,0.2\} & \{0.3,\,0.3,\,0.3\}\\
 \hline
 Average Travel Time [s] & 12.31 & 12.89 & 13.63 \\
 Average Time Gap [s] & 2.56 & 2.57 & 2.60\\
 Minimum Time Gap [s] & 1.26 & 1.52 & 1.73 \\
 Maximum Time Gap [s] & 4.56 & 4.56 & 4.56\\
\hline
\hline
\end{tabular*}
\end{table}

\begin{figure}[t!]
\centering
\includegraphics[width=.9\columnwidth]{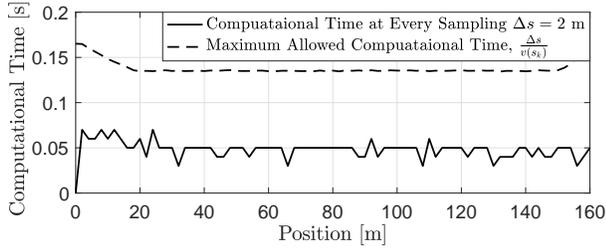}\\[-2ex]
\caption{Computational time of an example CAV $i$ and the corresponding maximum allowed computational time $\Delta s/v_i(k)$ with sampling $\Delta s\!=\!2$\,m by solving the decentralized HRCS problem.}
\label{fig:computaional_time}
\end{figure}

The proposed HRCS computational time of a single CAV $i$ at every step with the sampling interval $\Delta s\!=\!2$\,m
is shown in Fig.~\ref{fig:computaional_time}. 
The dashed line denotes the estimated permissible computational time at every interval distance $\Delta s$ in the { spatial domain}, obtained by $\Delta s/v_i(k)$. As can be seen, the computational time of vehicle $i$ is strictly below the maximum allowed time, which validates the implementation { potential} of the proposed HRCS approach.

\section{Conclusions}\label{sec:conclusions}
This paper proposes a new hierarchical robust control strategy (HRCS) for decentralized autonomous intersection coordination of connected and autonomous electric vehicles. The HRCS determines the optimal crossing order and velocity trajectories successively using an optimal control method and a tube-based robust model predictive control (RMPC) in a decentralized traffic coordination scheme. In particular, the RMPC can cope with the additive disturbances and modeling uncertainties entailed in the vehicle dynamic model and onboard sensor measurements, and therefore, it is the key to securing safety in reality. The optimization problems entailed in the control framework are solved as convex second-order cone programs with a suitable relaxation and spatial modeling approach, which can guarantee a fast and unique solution. The equivalence between the relaxed and the original problems is validated. 

Numerical examples verify the effectiveness and robustness of the proposed HRCS by comparisons against a nominal MPC-based strategy. The energy-time trade-off is examined for different prediction horizon lengths and arrival rates by finding the Pareto front of optimal solutions. The Pareto front shows an increase of 20\% in travel time can lead to an average energy consumption reduction of about 24\% at an arrival rate of 800 veh/h per lane with prediction horizon length $N_p\!=\!15$.
Finally, the computational efficiency of the convex HRCS is examined for a distance interval $\Delta s\!=\!2$~m and the results show the practical potential of the proposed scheme.



\bibliographystyle{IEEEtran}
\bibliography{paper}

\end{document}